\documentclass[12pt]{article}
\usepackage{graphicx}
\usepackage{psfrag}
\usepackage{mathrsfs}
\usepackage{amssymb}
\usepackage{amsfonts}
\usepackage{latexsym}
\usepackage{amsmath}
\usepackage{amsthm}
\usepackage[cp1251]{inputenc}
\usepackage[english]{babel}

\newcommand{\er}[1]{\textrm{(\ref{#1})}}
\def\lb{\label}

\theoremstyle{plain}
\newtheorem{Def}{\bf Definition}
\newtheorem{theorem}{\bf Theorem}[section]
\newtheorem{lemma}[theorem]{\bf Lemma}

\newtheorem{corollary}[theorem]{\bf Corollary}
\theoremstyle{remark}

\def\a{\alpha}                \def\A\mathbf{A}
\def\b{\beta}                 \def\B\mathbf{B}
\def\g{\gamma}           \def\cC{\mathcal{C}}     \def\bC\mathbf{C}
\def\G{\Gamma}                \def\bD\mathbf{D}
                \def\bE\mathbf{E}
\def\D{\Delta}                \def\bF\mathbf{F}
\def\ve{\varepsilon}     \def\cG{\mathcal{G}}     \def\bG\mathbf{G}
\def\z{\zeta}            \def\cH{\mathcal{H}}     \def\bH\mathbf{H}
\def\e{\eta}                  \def\bI\mathbf{I}
\def\vt{\vartheta}            \def\bJ\mathbf{J}
\def\vT{\Theta}          \def\cK{\mathcal{K}}     \def\bK\mathbf{K}
                \def\bL\mathbf{L}
\def\l{\lambda}          \def\cM{\mathcal{M}}     \def\bM\mathbf{M}
               \def\bN\mathbf{N}
\def\m{\mu}                   \def\bO\mathbf{O}
\def\n{\nu}              \def\cP{\mathcal{P}}     \def\bP\mathbf{P}
\def\r{\rho}                  \def\bQ\mathbf{Q}
\def\s{\sigma}           \def\cR{\mathcal{R}}     \def\bR\mathbf{R}
           \def\cS{\mathcal{S}}     \def\bS\mathbf{S}
\def\t{\tau}             \def\cT{\mathcal{T}}     \def\bT\mathbf{T}
\def\f{\phi}                  \def\bU\mathbf{U}
\def\F{\Phi}                  \def\bV\mathbf{V}
\def\vp{\varphi}              \def\bW\mathbf{W}
\def\c{\chi}                  \def\bX\mathbf{X}
\def\p{\psi}                  \def\bY\mathbf{Y}
             \def\cZ{\mathcal{Z}}     \def\bZ\mathbf{Z}
\def\o{\omega}
\def\O{\Omega}
\def\x{\xi}

\def\vk{\varkappa}

\def\mA{{\mathscr A}}
\def\mB{{\mathscr B}}

\def\mH{{\mathscr H}}

\def\mP{{\mathscr P}}

\newcommand{\gD}{\mathfrak{D}}

\newcommand{\gF}{\mathfrak{F}}

\def\J{\mathbb{J}}
\def\Z{\mathbb{Z}}
\def\R{\mathbb{R}}
\def\C{\mathbb{C}}

\def\N{\mathbb{N}}

\def\qqq{\qquad}
\def\qq{\quad}
\let\ge\geqslant
\let\le\leqslant
\let\geq\geqslant
\let\leq\leqslant
\newcommand{\ca}{\begin{cases}}
\newcommand{\ac}{\end{cases}}

\newcommand{\ma}{\begin{pmatrix}}
\newcommand{\am}{\end{pmatrix}}

\def\lt{\biggl}
\def\rt{\biggr}
\let\geq\geqslant
\let\leq\leqslant
\def\lra{\Leftrightarrow}
\def\[{\begin{equation}}
\def\]{\end{equation}}
\def\wt{\widetilde}
\def\wh{\widehat}
\def\pa{\partial}
\def\sm{\setminus}
\def\es{\emptyset}
\def\sps{\supset}
\def\ss{\subset}
\def\no{\noindent}
\def\ol{\overline}
\def\iy{\infty}
\def\ev{\equiv}
\def\/{\over}
\def\ts{\times}

\def\os{\oplus}

\def\Im{\mathop{\rm Im}\nolimits}
\def\supp{\mathop{\rm supp}\nolimits}

\def\Tr{\mathop{\rm Tr}\nolimits}

\def\BBox{\hspace{1mm}\vrule height6pt width5.5pt depth0pt \hspace{6pt}}

\begin{document}
\title{ Schr\"odinger operators on zigzag graphs}
\author{
Evgeny Korotyaev
\begin{footnote} {
Institut f\"ur  Mathematik,  Humboldt Universit\"at zu Berlin,
e-mail: evgeny@math.hu-berlin.de\ \ 
}
\end{footnote}
\and Igor Lobanov
\begin{footnote} {
Mathematical Faculty, Mordovian State University, 430000 Saransk, 
e-mail: lobanov@math.mrsu.ru }
\end{footnote}
}
\maketitle

\begin{abstract}
\no We consider the Schr\"odinger operator on zigzag  
 graphs with a periodic potential. The spectrum of this operator consists of an absolutely continuous part (intervals separated by gaps) plus an infinite number of eigenvalues  with infinite multiplicity. We describe all compactly supported eigenfunctions 
with the same eigenvalue. We define a Lyapunov function,
which is analytic on some Riemann surface. On each
sheet, the Lyapunov function has the same properties 
as in the scalar case, but it has 
branch points, which we call resonances. We prove that all  resonances are real. We determine the asymptotics of the periodic and anti-periodic spectrum and of the resonances at high energy. We show that there exist two types of
gaps: i) stable gaps, where the endpoints are periodic and
anti-periodic eigenvalues, ii) unstable (resonance) gaps, where the
endpoints are resonances (i.e., real branch points of the Lyapunov
function). 
We obtain the following results from the inverse spectral theory:
1) we describe all finite gap potentials, 2) the mapping: potential -- all eigenvalues is a real analytic isomorphism for some class of potentials.

\no We apply all these results to quasi-1D models of
zigzag single-well carbon nanotubes.
\end{abstract}

\section{Introduction  and main results}
\setcounter{equation}{0}

\begin{figure}\lb{fig1}
\centering
\noindent
(a){
\tiny
\psfrag{g001}[l][l]{$\Gamma_{0,0,2}$}
\psfrag{g002}[l][l]{$\Gamma_{0,0,3}$}
\psfrag{g003}[l][l]{$\Gamma_{0,0,1}$}
\psfrag{g011}[c][c]{$\Gamma_{0,1,2}$}
\psfrag{g012}[c][c]{$\Gamma_{0,1,3}$}
\psfrag{g013}[c][c]{$\Gamma_{0,1,1}$}
\psfrag{g021}[c][c]{$\Gamma_{0,2,2}$}
\psfrag{g022}[c][c]{$\Gamma_{0,2,3}$}
\psfrag{g023}[c][c]{$\Gamma_{0,2,1}$}
\psfrag{g-101}[c][c]{$\Gamma_{-1,0,2}$}
\psfrag{g-102}[l][l]{$\Gamma_{-1,0,3}$}
\psfrag{g-103}[l][l]{$\Gamma_{-1,0,1}$}
\psfrag{g-111}[c][c]{$\Gamma_{-1,1,2}$}
\psfrag{g-112}[c][c]{$\Gamma_{-1,1,3}$}
\psfrag{g-113}[c][c]{$\Gamma_{-1,1,1}$}
\psfrag{g-121}[c][c]{$\Gamma_{-1,2,2}$}
\psfrag{g-122}[c][c]{$\Gamma_{-1,2,3}$}
\psfrag{g-123}[c][c]{$\Gamma_{-1,2,1}$}
\psfrag{g101}[l][l]{$\Gamma_{1,0,2}$}
\psfrag{g102}[l][l]{$\Gamma_{1,0,3}$}
\psfrag{g103}[l][l]{$\Gamma_{1,0,1}$}
\psfrag{g111}[c][c]{$\Gamma_{1,1,2}$}
\psfrag{g112}[c][c]{$\Gamma_{1,1,3}$}
\psfrag{g113}[c][c]{$\Gamma_{1,1,1}$}
\psfrag{g121}[c][c]{$\Gamma_{1,2,2}$}
\psfrag{g122}[c][c]{$\Gamma_{1,2,3}$}
\psfrag{g123}[c][c]{$\Gamma_{1,2,1}$}
\includegraphics[width=.65\textwidth,height=.5\textwidth]{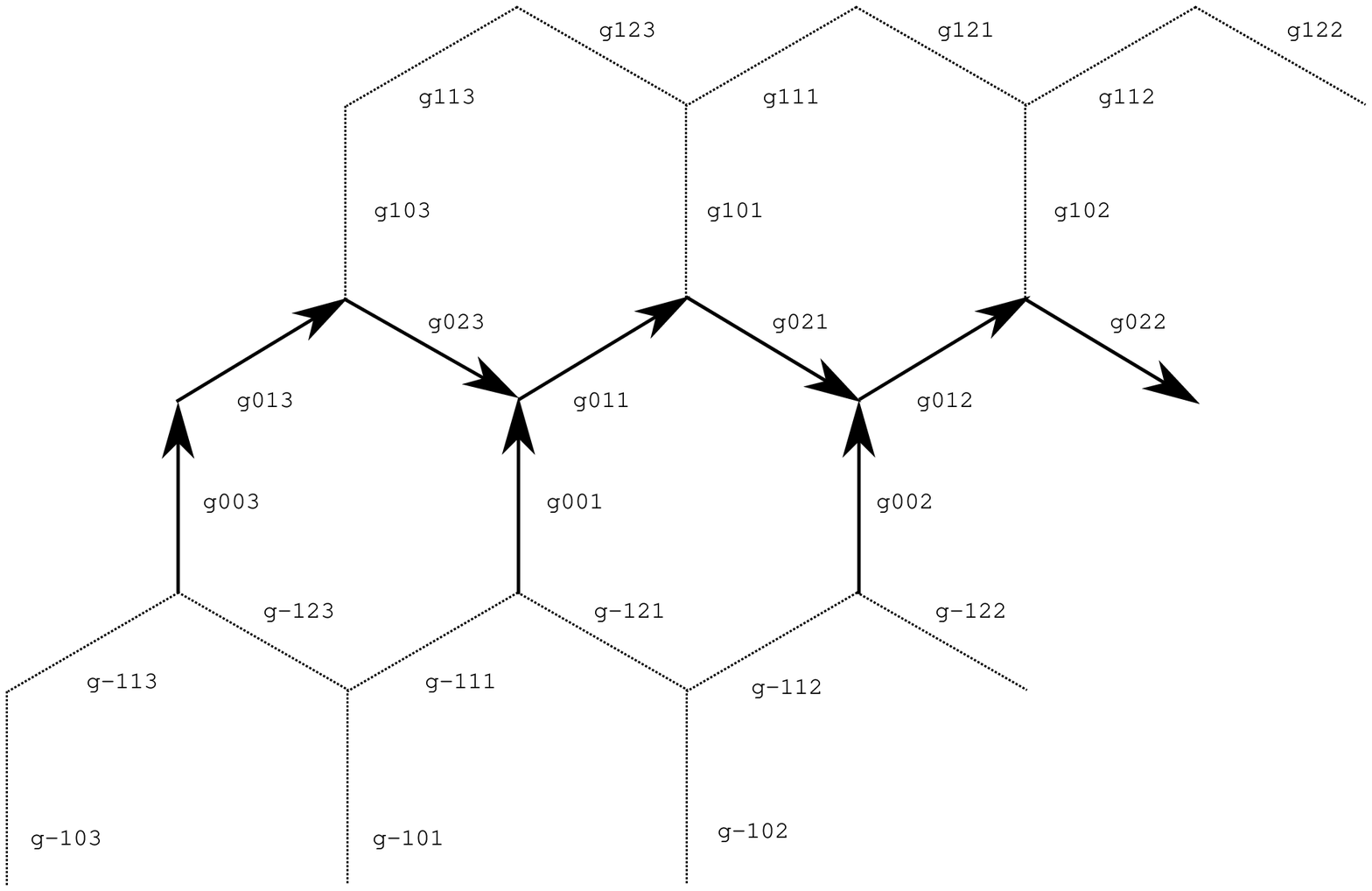}
}
(b){
\tiny
\psfrag{g-10}[l][l]{$\Gamma_{-1,0}$}
\psfrag{g-11}[l][l]{$\Gamma_{-1,1}$}
\psfrag{g-12}[l][l]{$\Gamma_{-1,2}$}
\psfrag{g00}[l][l]{$\Gamma_{0,0}$}
\psfrag{g01}[l][l]{$\Gamma_{0,1}$}
\psfrag{g02}[l][l]{$\Gamma_{0,2}$}
\psfrag{g10}[c][c]{$\Gamma_{1,0}$}
\psfrag{g11}[c][c]{$\Gamma_{1,2}$}
\psfrag{g12}[c][c]{$\Gamma_{1,3}$}
\includegraphics[height=.5\textwidth]{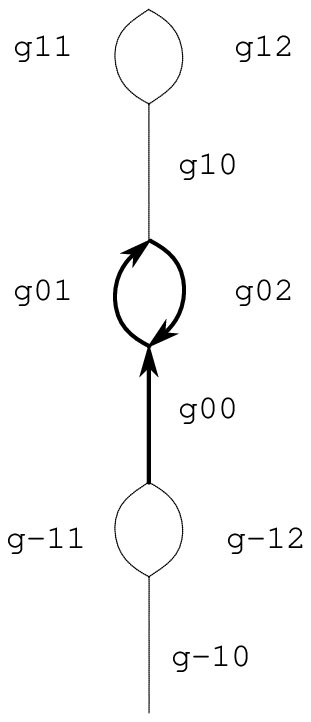}
}
\caption{Zigzag graph $\G^N$ for (a) $N=3$; (b) $N=1$.
The fundamental domain $\Gamma_0$ is marked by a bold line.}
\end{figure}

\begin{figure}
\noindent
\centering
 
\tiny
\psfrag{e0}{$e_0$}
\psfrag{e1}{$e_1$}
\psfrag{e2}{$e_2$}
\psfrag{0r01}{$r_{0,0,1}$}
\psfrag{1r01}{$r_{0,1,1}$}
\psfrag{0r02}{$r_{0,0,2}$}
\psfrag{0r03}{$r_{0,0,3}$}
\psfrag{0r04}{$r_{0,0,4}$}
\psfrag{0r0-1}{$r_{0,0,-1}$}
\psfrag{0r00}{$r_{0,0,0}$}
\psfrag{0r01}{$r_{0,0,1}$}
\psfrag{0r11}{$r_{1,0,1}$}
\psfrag{0r12}{$r_{1,0,2}$}
\psfrag{0r13}{$r_{1,0,3}$}
\psfrag{0r21}{$r_{2,0,1}$}
\psfrag{0r22}{$r_{2,0,2}$}
\psfrag{0r23}{$r_{2,0,3}$}
\psfrag{0r-11}{$r_{-1,0,1}$}
\psfrag{0r-12}{$r_{-1,0,2}$}
\psfrag{0r-13}{$r_{-1,0,3}$}
\psfrag{0r-21}{$r_{-2,0,1}$}
\psfrag{0r-22}{$r_{-2,0,2}$}
\psfrag{0r-23}{$r_{-2,0,3}$}
\includegraphics[width=.65\textwidth,height=.5\textwidth]{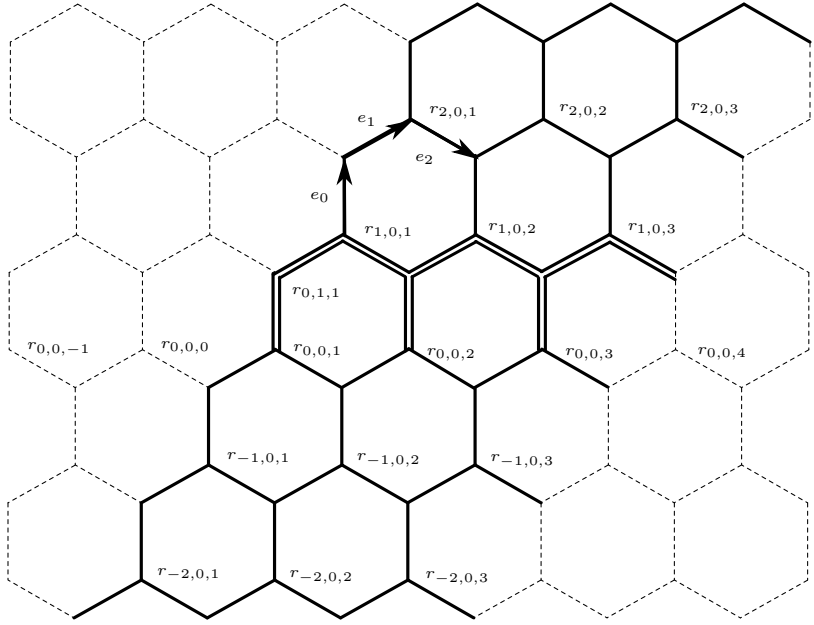}
\caption{The unrolled honeycomb lattice. The zigzag graph
is marked by bold, N=3. The fundamental domain is marked by a double bold line. Each edge $\G_\o$ has orientation 
given by $e_\o=e_j, j=0,1,2$ with the starting point $r_\o\in \R^2$.
}
\label{fig2}
\end{figure}

Consider the Schr\"odinger operator $H=-\pa^2+q$
with a periodic potential $q$ on the zigzag  graph $\G^N$.
The zigzag  graph $\G^N$ is related to a zigzag  nanotube $T^N$.
We shall discuss this at the end of this section.
In order  to define $\G^N$ we introduce the so-called honeycomb lattice $\G=\cup_{\o} \G_\o\ss\R^2$, $\o=(n,j,k)\in \Z\ts \J\ts \Z, \J=\{0,1,2\}$, where the edges $\G_\o$ are given by
$$
\G_\o=\{x=r_\o+te_\o,\qq  t\in [0,1]\},\qq  e_\o=e_j\in \R^2,\  
\qq r_{(n,1,k)}=r_{(n,0,k)}+e_0,
$$
$$
r_{(n,0,k)}=r_{(n-1,2,k)}=n(e_0+e_1)+k(e_1+e_2),\qq 
e_0=(0,1),\ e_1={1\/2}(\sqrt 3,1),\ e_2={1\/2}(\sqrt 3,-1),\ 
$$
see Fig. \ref{fig1} and \ref{fig2}. The edges $\G_\o$ of length 1 form regular hexagons. Each edge $\G_\o$ has orientation 
given by $e_\o=e_j$ with the starting point $r_\o\in \R^2$.
We have the coordinate $x=r_\o+te_\o$ and the local coordinate $t\in [0,1]$. We define an oriented  zigzag graph $\G^N$ by identifying $\G_{n,j,k}$ with $\G_{n,j,k+N}$ for all $(n,j,k)\in\Z\ts\J\ts \Z$.
Here and below for simplicity we write $\G_{n,j,k}=\G_{(n,j,k)}$ and $f_{n,j,k}=f_{(n,j,k)}$ etc. Thus $\G^N$ is a topological space in the quotient topology.   As example, the graphs $\G^3$ and $\G^1$ are illustrated in Fig. 1.  We represent $\G^N$ by 
$$
\G^N=\cup_{\o\in \cZ} \G_\o,\qq
\qq\o=(n,j,k)\in \cZ=\Z\ts \J\ts \Z_N,\qq \Z_N=\Z/(N\Z),\qq
\J=\{0,1,2\}. 
$$
In our paper we assume that $N=2m+1$ for some integer $m\ge 0$. 
For each function $y$ on $\G^N$ we define a function $y_\o=y|_{\G_\o}, \o\in \cZ$. We identify each function $y_\o$ on $\G_\o$ with a function on $[0,1]$ by using the local coordinate $t\in [0,1]$ and define the Hilbert space $L^2(\G^N)=\os_{\o\in \cZ} L^2(\G_\o)$. 
Introduce the space $C(\G^N)$  of continuous functions on $\G^N$
 and the Sobolev space $W^2(\G^N)=\{y, \os_{\o\in \cZ} y_\o''\in L^2(\G^N); \ y $ {\it satisfies  
{\bf the  Kirchhoff Boundary Conditions:} $y\in C(\G^N)$ and 
 the following identities hold
\[
\lb{KirC}
-y_{n,0,k}'(1)+y_{n,1,k}'(0)-y_{n,2,k-1}'(1)=0,\qqq 
y_{n+1,0,k}'(0)-y_{n,1,k}'(1)+y_{n,2,k}'(0)=0.
\]
for all $(n,k)\in \Z\ts\Z_N$.}\} 
Condition \er{KirC} means  that  the sum of derivatives of $y$ at each vertex  of $\G^N$ equals 0 and the orientation of edges gives
the sign $\pm$. Our operator $H$ on $\G^N$ acts in the Hilbert space $L^2(\G^N)$. We define our operator $H$ by
$(Hy)_\o=-y_\o''+q y_\o$, where $\ y=(y_\o)_{\o\in \cZ}\in \gD(H)=W^2(\G^N),\ q\in L^2(0,1).
$
In the case $q=0$, we denote the operator $H$ by $H^0$.
The operator $H^0$ is self-adjoint, see [Ca1], and so is $H$,
 see Sect.3.

Recall the needed properties of the Hill operator $\wt H y=-y''+q(x)y$ on 
the real line with a periodic potential $q(x+1)=q(x),x\in \R$.
The spectrum of $\wt H$ is purely absolutely continuous and
consists of intervals $\wt\s_n=[\l_{n-1}^+,\l_n^-], n\ge 1$. These intervals are separated by the gaps $\g_n=(\l_n^-,\l_n^+)$ of length $|\g_n|\ge 0$. If a gap $\g_n$ is degenerate, i.e. $|\g_n|=0$, then the corresponding segments $\wt\s_n,\wt\s_{n+1}$ merge. For the equation $-y''+q(x)y=\l y$ on the real line we define the fundamental solutions $\vt(x,\l)$ and $\vp(x,\l),x\in \R$ satisfying $\vt(0,\l)=\vp'(0,\l)=1, \vt'(0,\l)=\vp(0,\l)=0$. The corresponding monodromy matrix $\cM$
and the Lyapunov function $\D$  are given by
\[
\cM(\l)=\ma\vt(1,\l) & \vp(1,\l) \\ 
\vt'(1,\l) & \vp'(1,\l)\am,\ \ \ \ \ 
\D(\l)={\vp'(1,\l)+\vt(1,\l)\/2},\ \ \l\in\C.
\]
The sequence $\l_0^+<\l_1^-\le \l_1^+\ <.....$
is the spectrum of the equation $-y''+qy=\l y$ with periodic
boundary conditions of period 2,  that is  $y(x+2)=y(x), x\in \R$.
 Here equality $\l_n^-= \l_n^+$ means that $\l_n^\pm$ is an eigenvalue of multiplicity 2. Note that $\D(\l_{n}^{\pm})=(-1)^n, \  n\ge 1$. The lowest  eigenvalue $\l_0^+$ is simple, $\D(\l_0^+)=1$, and the
corresponding eigenfunction has period 1. The eigenfunctions
corresponding to $\l_n^{\pm}$ have period 1 if $n$ is even,
and they are anti-periodic, that is $y(x+1)=-y(x),\ x\in \R$, if
$n$ is odd. The derivative of the Lyapunov function has a zero $\l_n$ in each ''closed gap'' $[\l^-_n,\l^+_n]$, that is $ \D'(\l_n)=0$.  Let $\m_n, n\ge 1,$ be the spectrum of the problem $-y''+qy=\l y, y(0)=y(1)=0$ (the Dirichlet spectrum), and let $\n_n, n\geq 0,$ be the spectrum of the problem $-y''+qy=\l y, y'(0)=y'(1)=0$ (the Neumann  spectrum). It is well-known that $\m_n, \n_n \in [\l^-_n,\l^+_n ]$ and $\n_0\le \l^+_0$.
Moreover, a potential $q$ is even, i.e., $q\in L^2_{even}(0,1)=\rt\{q\in L^2(0,1): q(1-x)=q(x),x\in[0,1]\rt\}$ iff $|\g_n|=|\m_n-\n_n| $ for all $n\ge 1$. Define the set $\s_D=\{\m_n, n\ge 1\}$ and note that $\s_D=\{\l\in\C: \vp(1,\l)=0\}$.

For simplicity we shall denote $\G_{\a,1}\ss \G^1$ by  $\G_{\a}$, for $\a=(n,j)\in \cZ_1=\Z\ts \J$. Thus $\G^1=\cup_{\a\in \cZ_1} \G_\a$,
see Fig \ref{fig1}.
In Theorem \ref{T1} we will show that 
 $H$ is unitarily equivalent to $\os_1^N H_k$, where the operator $H_k$ acts in the Hilbert space $L^2(\G^1)$ and is given by 
\[
\lb{Hk}
(H_k f)_\a=-f_\a''+q f_\a,\qqq \os_{\a\in \cZ_1} f_\a''\in L^2(\G^1),
\]
and the vector function $f=(f_\a)_{\a\in\cZ_1}$ satisfies the Kirchhoff conditions:
\[
\lb{1K0}
f_{n,0}(1)=f_{n,1}(0)=s^k f_{n,2}(1),\quad
f_{n+1,0}(0)=f_{n,1}(1)=f_{n,2}(0),\qq s=e^{i{2\pi \/N}},
\]
\[
\lb{1K1}
-f'_{n,0}(1)+f'_{n,1}(0)-s^k f'_{n,2}(1)=0,\qq
f'_{n+1,0}(0)-f'_{n,1}(1)+f'_{n,2}(0)=0.
\]
We reduce the spectral problem on the graph $\G^1$ to some matrix problem on $\R$. In order to describe this we define the fundamental subgraph $\G_0^1$ of $\G^1$ by $\G_{0}^1=\cup_{j=0}^2 \G_{0,j}$, see Fig. \ref{fig1}. On $\G^1$, the group $\Z$ acts via
$
p\circ  \G_{n,j}^1=\G_{n+p,j}^1,\  (n,j,p)\in \cZ_1\ts\Z. 
$
Thus $\G_0^1$ is a fundamental domain associated with this group action
of $\Z$.

For the operator $H_k$ we construct the fundamental solutions 
$
\vT_k(x,\l)
=(\vT_{k,\a}(x,\l))_{\a\in\cZ_1},$ and $ \F_k(x,\l)=(\F_{k,\a}(x,\l))_{\a\in\cZ_1}$, $(x,\l)\in \R\ts \C,
$ 
which  satisfy
\[
\lb{eqf}
-f_\a''+qf_\a=\l f_\a,\ \ 
\ \ \text{the Kirchhoff Boundary Conditions \er{1K0},\er{1K1}}, 
\]
\[
\lb{eqf0}
\vT_{k,\b}(0,\l)=\F_{k,\b}'(0,\l)=1,\qq 
\vT_{k,\b}'(0,\l)=\F_{k,\b}(0,\l)=0,\qqq \b=(0,0).
\]
We introduce the monodromy matrix
\[
\cM_k(\l)=\ma \vT_{k,\a}(0,\l) & \F_{k,\a}(0,\l)\\
\vT_{k,\a}'(0,\l) & \F_{k,\a}'(0,\l)\am,\qqq \a=(1,0),
\]
which is determined by $\vT_k, \F_k$ on the fundamental domain
$\G_0^1$.

There are two methods to study periodic differential operators.
The direct integral analysis usually used for partial
differential operators [ReS] gives general information about the spectrum, but no detailed results.
The method of ordinary differential operators, based on the Floquet matrix analysis, gives detailed results. Note that there are a lot of open problems [BBK] even for the Schr\"odinger operator with periodic $2\ts 2$ matrix potentials on the real line .

We introduce the monodromy matrix, similar to the case
of Schr\"odinger operators with periodic matrix potentials
on the real line, see [YS]. After this, using the approach from [BK],[BBK],[CK] we introduce the Laypunov functions and study   
the properties of this functions, similar to the 
case of  the Schr\"odinger operator with a periodic matrix
potential on the real line. This is a crucial point of our analysis.
 Here we essentially use the results and techniques
from the papers [BK], [BBK],[CK]. The recent papers [BBK],[CK] are devoted to Schr\"odinger operators with periodic matrix
potentials (a standard case) on the real line. Remark that Carlson \cite{Ca} studied the monodromy operator to analyze the Schr\"odinger operator  on  a product of graphs. His results do not cover 
our case.

We  formulate our first preliminary results.

\begin{theorem}\label{T1}
(i) The operator $H$ is unitarily equivalent to $\os_1^N H_k$, where the operator $H_k$ is given by \er{Hk}-\er{1K1}.

\no (ii) For any $\l\in\C\sm\s_D$ and $k\in \Z_N$  there exist unique fundamental solutions $\vT_k, \F_k$ of the system \er{eqf},\er{eqf0}.
 Moreover, each of the functions $\vT_k(x,\l),\F_k(x,\l),x\in \G^1$ is meromorphic in $\l\in\C\sm\s_D$ and each matrix $\cM_k(\l)$ satisfies
\[
\lb{T1-1}
\cM_k=\cR^{-1}\cT_k \cR \cM ,\quad
\cT_k={s^{-{k\/2}}\/2c_k}\ma 2\D & 1 \\ 4\D^2-4c_k^2 & 2\D \am,\qq 
\cR=\ma 1 & 0\\0 & \vp(1,\cdot)\am.
\]
In particular,
\[
\label{T1-2}
\det \cM_k={s}^{-k},\qq
\Tr\cM_k={2(\D_0+s_k^2)\/s^{{k\/2}}c_k},\qq
\D_0\ev {\Tr \cM_0\/2}=2\D^2+{\vt'(1,\cdot)\vp(1,\cdot)\/4}-1, 
\]
where $c_k=\cos{\pi k\/N},\ s_k=\sin{\pi k\/N} s=e^{i{2\pi\/N}}$.
Furthermore, $\cR\cM_k\cR^{-1}=\cT_k\cR\cM\cR^{-1}$ is an entire matrix-valued function. In particular, the function $D_k(\t,\l)=\det(\cM_k(\l)-\t I_{2})$, where $I_2$ is the $2\ts 2$  identity matrix, is entire with respect to $\l, \t\in \C$.

\end{theorem}

Remark that in contrast to the Schr\"odinger operator with periodic matrix potential on the real line (see \cite{YS} or \cite{CK}),
the monodromy matrix $\cM_k$ has poles at the points $\l\in \s_D$,
which are eigenvalues of $H_k$, see Theorem \ref{T2}.
Such a phenomenon has already been observed e.g. in 
\cite{SA}, \cite{MV}, \cite{Ku}.

Define the subspace $\cH_k(\l)=\{\p\in \gD(H_k): H_k\p=\l \p\}$
for $\l\in \s_{pp}(H_k), k\in \Z_N$. If 
$\dim \cH_k(\l_0)=\iy$ for some $\l_0\in \s_{pp}(H_k)$, then we say that $\{\l_0\}$ is a flat band. In Theorem \ref{T2} we describe all flat bands
and supports of eigenfunctions (see Fig. \ref{fig3}).

\begin{theorem}\lb{T2} 
Let $(\l,k)\in\s_D\ts\Z_N$. Then

\no (i) Each eigenfunction from $\cH_k(\l)$ vanishes at all vertices of $\G^1$.

\no (ii) Let $\vp=\vp(\cdot,\l)|_{[0,1]}$ and put $c=\vp'(1,\l)$. Define the function $\p^{(0)}=(\p^{(0)}_\a)_{\a\in \cZ_1}$  by :

if $\e=1-s^kc^2\ne 0$, then
\begin{multline}
\lb{T2-1}
\p^{(0)}_{n,j}=0, \text{for all}\  n\ne 0,-1, \ j\in \Z_3,\qqq
and \qq
\p^{(0)}_{0,0}=\e \vp, \ \p^{(0)}_{0,1}=c\vp, \ 
\p^{(0)}_{0,2}=c^2\vp, \\
\p^{(0)}_{-1,0}=0,\ 
\p^{(0)}_{-1,1}=-s^kc\vp, \ \p^{(0)}_{-1,2}=-\vp,
\end{multline}

if $\e=0$, then
\[
\lb{T2-2}
\p^{(0)}_{0,0}=0,\qq \p^{(0)}_{0,1}=\vp,\qq
\p^{(0)}_{0,2}=c\vp,\qq 
\p^{(0)}_{n,j}=0,\ all \ n\ne 0, j\in \Z_3.
\]
Then  each $\p^{(n)}=(\p^{(0)}_{n-m,j})_{(m,j)\in \cZ_1}\in \cH_k(\l), n\in \Z$  and each $f\in \cH_k(\l)$ has the form
\[
\lb{T2-3}
f=\sum_{n\in\Z}\wh f_n\p^{(n)},\qqq 
\wh f_{n}=\ca\e^{-1} f_{n,0}'(0) & if \ \e\ne 0\\
f_{n,1}'(0) & if \ \e= 0 \ac , 
\qqq (\wh f_n)_{n\in\Z}\in \ell^2.
\]
Moreover, the mapping $f\to \{\wh f_{n}\}_{n\in\Z}$ is a linear isomorphism between $\cH_k(\l)$ and $\ell^2$.

\end{theorem}

\begin{figure}
\tiny
\centering
\noindent
(a)\quad\quad\includegraphics[angle=90,width=.4\textwidth]{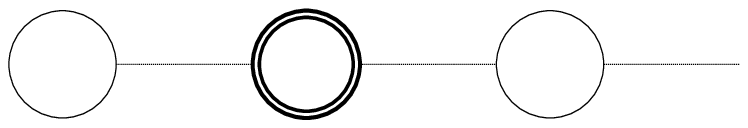}
(b)\quad\quad\includegraphics[angle=90,width=.4\textwidth]{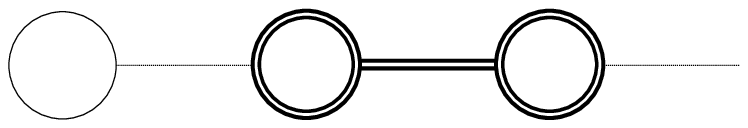}
\caption{The supports of the eigenfunction $\p^{(0)}$: (a) $\e=1-s^k\vp'(1,\l)^2=0$; (b) $\e\neq 0$.
}
\lb{fig3}
\end{figure}

\begin{figure}
\centering
\tiny
\psfrag{-5/4}[r]{$-\frac{5}{4}$}
\psfrag{-1}[r]{$-1$}
\psfrag{1}[r]{$1$}
\psfrag{0}[r]{$0$}
\noindent
\begin{tabular}{cc}
{
\psfrag{l0+}{$\l_{0,0}^+$}
\psfrag{l1-}{$\l_{0,1}^-$}
\psfrag{l1+}{$\l_{\!0,\!1}^+$}
\psfrag{l2-}{$\l_{\!\!0,\!2}^-$}
\psfrag{l2+}{$\l_{0,\!2}^+$}
\psfrag{l3-}{$\l_{0,\!3}^-$}
\psfrag{l3+}{$\l_{0,\!3}^+$}
\psfrag{d1}{$\l_{0,\!1}$}
\psfrag{d2}{}
\psfrag{d3}{}
\psfrag{z1}{$\eta_{0,\!1}$}
\psfrag{z2}{$\eta_{0,\!2}$}
\psfrag{z3}{$\eta_{0,\!3}$}
\psfrag{D0}{$\Delta_0$}
\includegraphics[width=.45\textwidth]{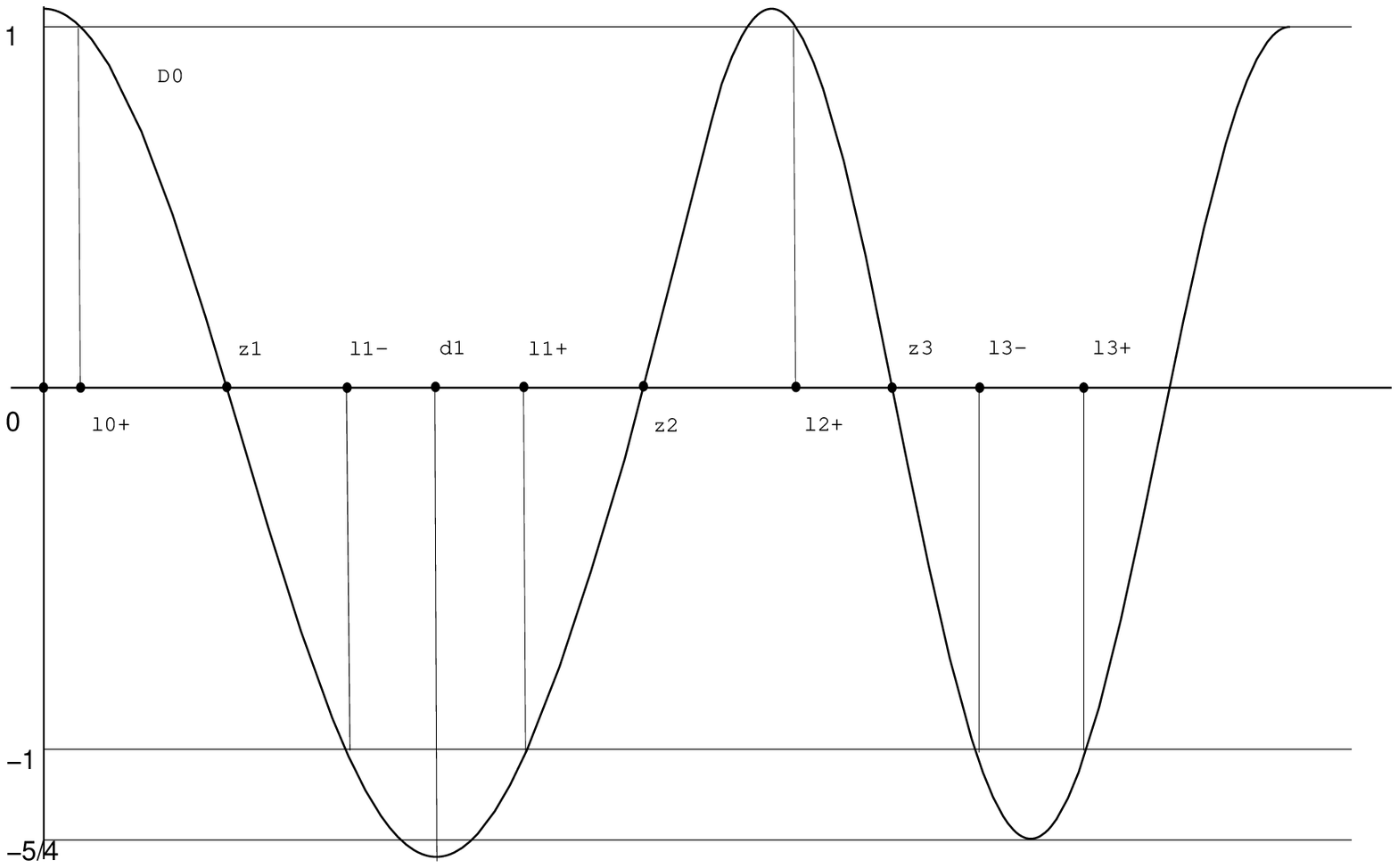}
}&
{
\psfrag{r0+}[c]{$r_{1,0}^+$}
\psfrag{r1-}[c]{$r_{1,1}^-$}
\psfrag{r1+}[c]{$r_{1,1}^+$}
\psfrag{r2-}[c]{$r_{1,2}^-$}
\psfrag{r2+}[c]{$r_{1,2}^+$}
\psfrag{r3-}[c]{$r_{1,3}^-$}
\psfrag{r3+}[c]{$r_{1,3}^+$}
\psfrag{r4-}[c]{$r_{1,4}^-$}
\psfrag{d1}{}
\psfrag{d2}{}
\psfrag{d3}{}
\psfrag{D1}{$\Delta_1$}
\includegraphics[width=.45\textwidth]{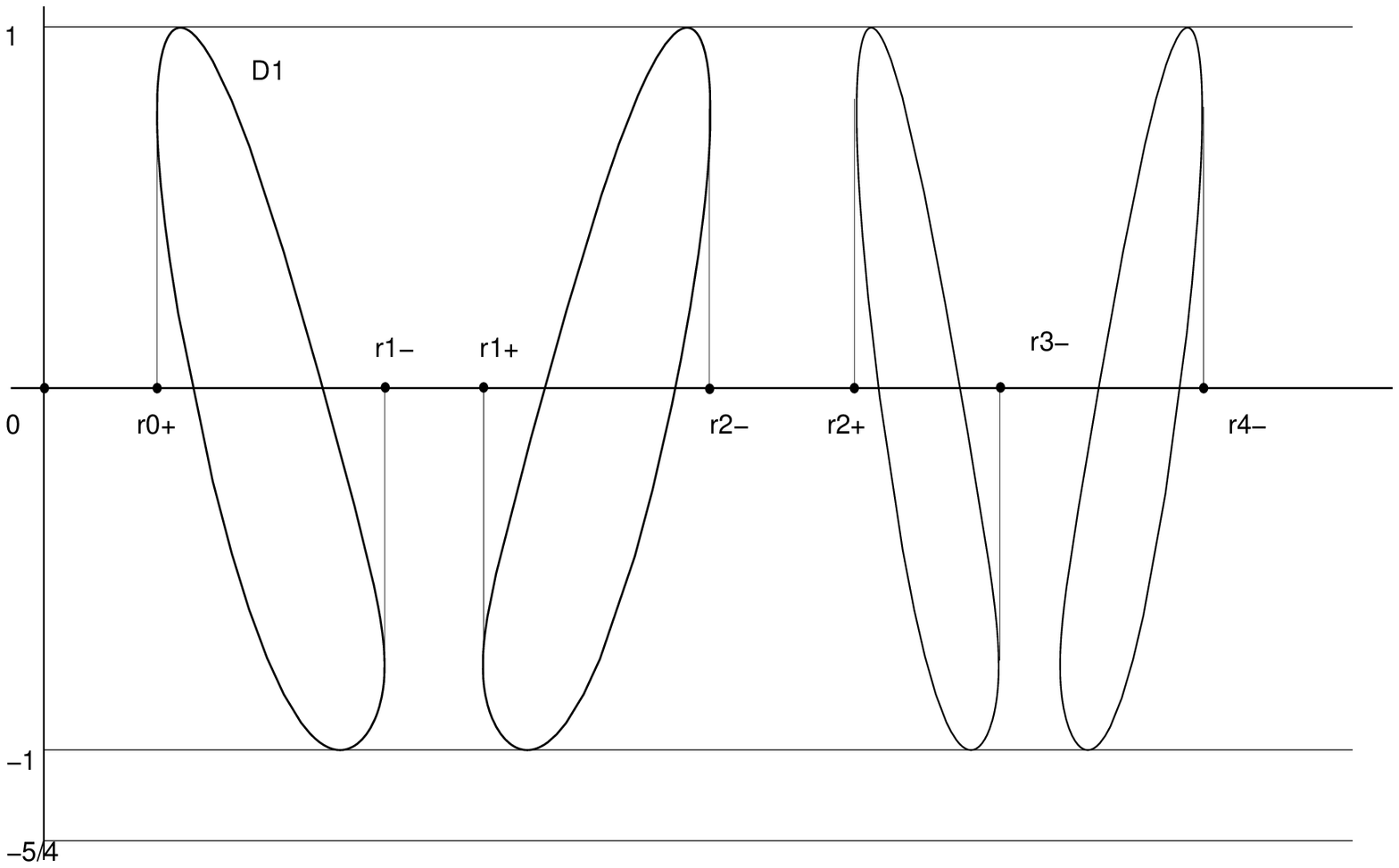}
}\\
{
\psfrag{r0+}[c]{$r_{3,0}^+$}
\psfrag{r2-}{$r_{3,2}^-$}
\psfrag{r2+}[c]{$r_{3,2}^+$}
\psfrag{r4-}{$r_{3,4}^-$}
\psfrag{d1}{}
\psfrag{d2}{}
\psfrag{d3}{}
\psfrag{D3}{$\Delta_3$}
\includegraphics[width=.45\textwidth]{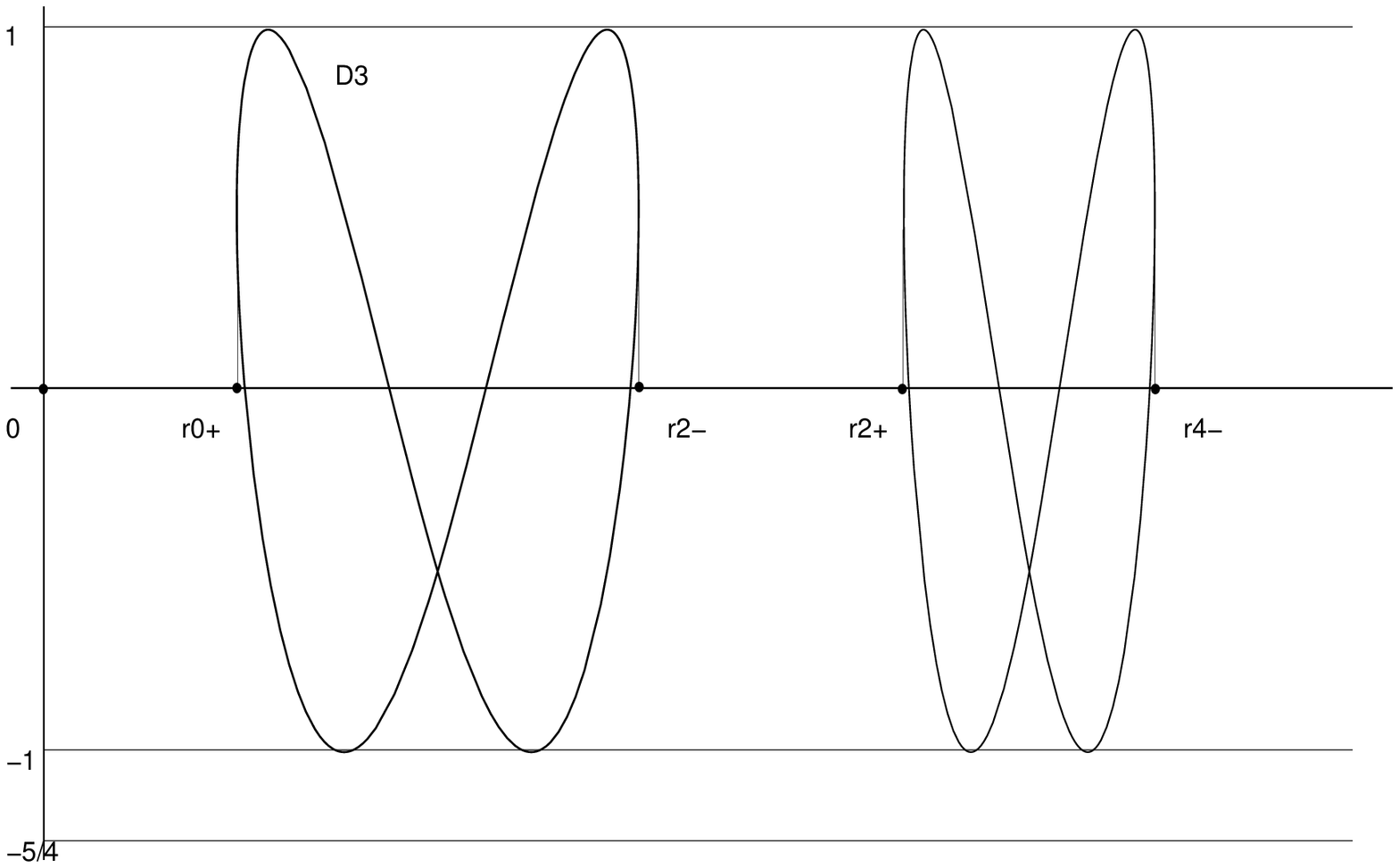}
}&
\includegraphics[width=.3\textwidth]{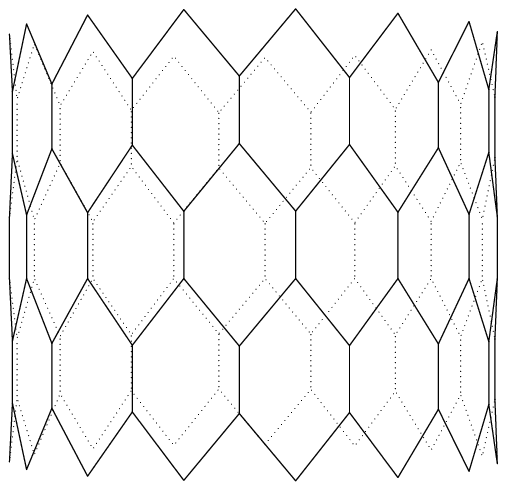}
\end{tabular}
\caption{Typical Lyapunov functions and the zigzag nanotube.
The function $\Delta_3$ is shown for $q\in L_\mathrm{even}(0,1), N=6$. 
}
\label{fig4}
\end{figure}

Let $\t_{k,\pm}$ be the eigenvalues of $\cM_k, k\in \Z_N$. Using \er{T1-2} we deduce that
\[
\lb{malk}
\t_{k,-}\t_{k,+}=s^{-k}, \qq \t_{k,-}+\t_{k,+}=\Tr 
\cM_{k}=s^{-k}\Tr\cM_{-k},
\qq \t_{-k,+}=\t_{k,-}s^{k},\qq \t_{-k,-}=\t_{k,+}s^{k}.
\]
If $k=0$, then we introduce the standard Lyapunov function $\D_0={\Tr \cM_0\/2}={1\/2}(\t_{0,+}+{1\/\t_{0,+}})$,
which is entire. If $k\neq 0$, then we define the Lyapunov functions 
$\D_{k,\pm}={1\/2}(\t_{k,\pm}+{1\/\t_{k,\pm}})$ (see [BBK],[CK]) and
using  \er{malk} we get $\D_{-k,\pm}=\D_{k,\pm}$.
Below we prove the following identities
\[
\lb{DeL1}
\D_{k,\pm}
=\x_k\pm\sqrt{\rho_k},\qqq \x_k=\D_0+s_k^2,\qq \r_k={s_k^2\/c_k^2}(c_k^2-\x_k^2),\qq c_k=\cos{\pi k\/N},\ s_k=\sin{\pi k\/N}.
\]
 If $q=0$, then we denote the corresponding functions by
$\D_k^0, \r_k^0,..$. In particular, we have
\[
\D_0^0={9\cos 2\sqrt{\l}-1\/8},\ \qqq  \x_k^0=\D_0^0+s_k^2,\ \ \ \r_k^0={s_k^2\/c_k^2}(c_k^2-(\x_k^0)^2),
\]
Introduce the two sheeted Riemann surface $\cR_k$
(of infinite genus) defined by $\sqrt {\r_k}$.
The functions $\D_{k,\pm},k\in \ol m=\{1,2,..,m\}$, 
are the branches of $\D_k=\x_k+\sqrt{\r_k}$  (see Fig.\ref{fig4})
on the Riemann surface $\cR_k$. 
We describe spectral properties of $H_k$ in terms of  $\D_k$.

\begin{theorem}
\lb{T3} 
(i) The Lyapunov functions $\D_{k,\pm}, k\in\ol m$ satisfy \er{DeL1}. 
 
\no (ii) For any $k\in\Z_N$ the  following identities hold:
\[
\lb{T3-1}
\s(H_k)=\s_{pp}(H_k)\cup\s_{ac}(H_k),\qq
\s_{pp}(H_k)=\s_D,\qq \s_{ac}(H_k)=\{\l\in\R: \D_k(\l)\in [-1,1]\}.
\] 
(iii) If $\D_k(\l)\in (-1,1)$ for some $\l\in \R,k=0,..,m$ and 
$\l$ is not a branch point of $\D_k(\l)$, then $\D_k'(\l)\neq 0$. 

\end{theorem}

{\bf Remark.} 1) If we know $\D_0$, then we determine all
$\D_k, \r_k,k\in\Z_N$ by \er{DeL1}. 2) If we know $\r_j$ for some $j\in\ol m$, then  by \er{DeL1}, we determine all $\D_k, \r_k,k\in\Z_N$.

\begin{Def}
A zero of $\r_k, k\in \ol m=\{1,2,..,m\}$ is called a {\bf resonance} of $H$. \end{Def}

Let $\l_{0,2n}^\pm$ and $\l_{0,2n+1}^\pm,n\ge 0$ be the zeros of $\det (\cM_0-I_2)$ and $\det (\cM_0+I_2)$ respectively (counted
with multiplicity). In Theorem \ref{T4}  we will show that the periodic  eigenvalues $\l_{0,2n}^\pm$ and the anti-periodic eigenvalues $\l_{0,2n+1}^\pm$ satisfy 
\[
\lb{epa0}
\D_0(\l_{0,n}^\pm)=(-1)^n,\ \ \ 
\l_{0,0}^+<\l_{0,1}^-< \l_{0,1}^+<\l_{0,2}^-\le \l_{0,2}^+<\l_{0,3}^-< \l_{0,3}^+<\l_{0,4}^-\le \l_{0,4}^+<....
\]
For $H_0$ we introduce the spectral bands $\s_{0,n}=[\l_{0,n-1}^+,\l_{0,n}^-]$ and 
the gaps $\g_{0,n}=(\l_{0,n}^-,\l_{0,n}^+)$, $n\ge 1$. Using the asymptotics of the fundamental solutions $\vt, \vp,$ we
obtain
\[ 
\lb{asD0}
\D_0(\l)=\D_0^0(\l)+{9q_0\/8}{\sin 2\sqrt{\l}\/\sqrt{\l}}+{O(e^{|\Im\sqrt{\l}|})\/|\l|},\ \qq q_0=\int_0^1 q(t)dt
\]
as $|\l|\to \iy$.
The following theorem describes the basic properties of $H_0$.

\begin{theorem}\label{T4}
The function $\D_0$ is entire and has the following properties:

\no (i) The function $\D_0'$ has only real simple zeros
$\l_{0,n},n\ge 1$,
which are separated by the simple zeros $\e_{0,n}$ of $\D_0$:
$\e_{0,1}<\l_{0,1}<\e_{0,2}<\l_{0,2}<\e_{0,3}<...$ and satisfy
\[
\lb{T4-1}
\n_n,\m_n\in \g_{0,2n},\ \ 
\D_0(\l_{0,2n})\ge 1,\ \ \ \   \  \ \ \ \ \ \ \   \D_0(\l_{0,2n-1})\le -{5\/4},\ \ \ \ \  for \ any \ n\ge 1.
\]
\no (ii) The periodic and anti-periodic eigenvalues $\l_{0,n}^\pm,n\ge 0$ satisfy \er{epa0} and have asymptotics
\[
\lb{T4-2}
\l_{0,1+2n}^\pm=(\pi{2n+1\/2}\pm \f)^2+q_0+o(n^{-1}),\ \ \qq
\  \l_{2n}^\pm=(\pi n)^2+q_0\pm \rt|{|\hat q_n|^2}-{\hat q_{sn}^2\/9}\rt|^{1\/2}+{O(1)\/n},
\]
where  $\f=\arcsin {1\/3}\in [0,{\pi\/2}]$
and $\hat q_n=\int_0^1q(t)e^{i2\pi nt}dt, \hat q_{sn}=\Im \hat q_{n}$. Moreover, 
 $(-1)^n\D_0(\l)\ge 1$ for all $\l\in [\l^-_{0,n},\l^+_{0,n}]$ and 
 $\l_{0,n}\in [\l^-_{0,n},\l^+_{0,n}]$.
 
 \no (iii) $|\g_n|=|\m_n-\n_n|$  iff $\g_n=\g_{0,2n}$.

\no (iv) $\g_{0,2n}\ss \g_n$ for all $n\ge 1$. Moreover,
for fixed $n\ge 1$ we have $\g_{0,2n}=\es$ iff $\g_{n}=\es$.

\no (v) $\D_0(\l_{1+2n})=-{5\/4}$ for all $n\ge 0$ iff 
$q\in L_{even}^2(0,1)$.

\end{theorem}

Remark that the last theorem gives the asymptotics of the gap length.

Let $\l_{k,2n}^\pm$ and $\l_{k,2n+1}^\pm,n\ge 0$ be the zeros of $\det (\cM_k-I_2)$ and $\det (\cM_k+I_2)$. 
Below we will show that the periodic  eigenvalues $\l_{k,2n}^\pm$ and the anti-periodic eigenvalues $\l_{k,2n+1}^\pm$ satisfy 
the equations $\D_k(\l_{k,n}^\pm)=(-1)^n$. In Theorem \ref{T5} we show that $\l_{k,n}^\pm$ satisfy  the equation
\[
\lb{eqpk}
\D_0(\l_{k,2n}^\pm)=\cos {2\pi k\/N},\ \qq \D_0(\l_{k,2n+1}^\pm)=
-1,\ \ \ \ \ \ k\in \ol m,
\]
and the labeling is given by: each $\l_{k,n}^\pm$ is double and
\[
\lb{espk}
\l_{0,0}^+<\l_{1,0}^+<\l_{2,0}^+<...<\l_{m,0}^+<\l_{0,1}^-
<\l_{0,1}^+
<\l_{m,2}^-<\l_{m-1,2}^-<...<\l_{0,2}^-\le \l_{0,2}^+<...,
\]
\[
\lb{espkI}
\l_{0,2n-1}^\pm=\l_{k,2n-1}^\pm,\ \ (k,n)\in \ol  m\ts\N. 
\]
The periodic eigenvalues $\l_{k,2n}^{\pm}, k\in \ol m=\{1,..,m\}, n\ge 0$ (i.e., $\D_k(\l_{k,2n}^{\pm})=1$ ) satisfy
\[
\lb{Tas-2}
\l_{k,2n}^{\pm}=(\pi n\pm \f_k)^2+q_0+{o(1)\/n},\ \ \qq 
\f_k={1\/2}\arccos {1+8c_{2k}\/9}\in [0,{\pi\/2}] \qq
as \ \ \ \ n\to \iy.
\]

Let $r_{k,n}^\pm, k\in \ol m, n\ge 0$ be the zeros of $\r_k$.
Below we will show that $r_{k,n}^\pm$ satisfy the equation
\[
\lb{eqr}
\D_0(r_{k,2n}^\pm )=c_k-s_k^2, \  \qq 
\D_0(r_{k,2n+1}^\pm )=-c_k-s_k^2, \ \qq
\ \ \ k\in \ol m,
\]
they are real and the labeling is given by
\[
\lb{esr}
r_{k,0}^+<
r_{k,1}^-<r_{k,1}^+< r_{k,2}^-<r_{k,2}^+ <..,\ \ \ \ k\neq {N\/3},
\]\[
\lb{esrI}
r_{k,0}^+<
r_{k,1}^-\le r_{k,1}^+< r_{k,2}^-<r_{k,2}^+ <r_{k,3}^-\le r_{k,3}^+..,\ \ \ \ k= {N\/3}.
\]
The resonances $r_{k,2n+s}^{\pm}, s=0,1$ (for $k\in \ol m,k\ne {N\/3}$) satisfy 
\[
\lb{Tas-4}
r_{k,n}^\pm=({\pi n\/2}\pm b_{k,n})^2+q_0+{o(1)\/n},\qqq b_{k,2n}=\f_{k,0},\qq b_{k,2n+1}={\pi\/2}-\f_{k,1}, 
\]
$\f_{k,s}={1\/2}\arccos {1+(-1)^s8c_{k}-8s_k^2\/9}\in [0,{\pi\/2}]$
 as $n\to \iy$. We describe the spectral properties of $H$.

\begin{theorem}\lb{T5}
 (i) The periodic and anti-periodic eigenvalues $\l_{k,n}^\pm, k\in \ol m=\{1,..,m\},n\ge 0$ satisfy Eq. \er{eqpk} and the relations \er{espk}-
 \er{Tas-2}.

\no (ii) The resonances  $r_{k,n}^\pm,k\in \ol m,n\ge 0$ satisfy Eq. \er{eqr} and estimates \er{esr}-\er{Tas-4}.

\no (iii) For $ n\ge 1, k\in \ol m$ the  following identities are fulfilled:
\[
\lb{T5-1}
\s_{ac}(H)=\cup_{n\ge 1} S_n=\cup_0^N \s_{ac}(H_k),
\ \  S_n=[E_{n-1}^+,E_n^-]=\cup_{k=0}^m \s_{k,n},\ \  \ 
\cap_{k=0}^m \s_{k,n}\ne \es,
\]
\[
\lb{T5-2}
\s_{0,n}=[\l_{0,n-1}^+,\l_{0,n}^-],\ \ \  {\rm and }\ \ \ 
\l_{k,n}^\pm\in \s_{k,n}=[r_{k,n-1}^+r_{k,n}^-],\ \ \ 
\]
\[
\lb{T5-3}
 G_n=(E_{n}^-,E_n^+)=\cap_0^m \g_{k,n},\ \ G_{2n}=\g_{0,2n}.
\]
\no (iv) If  $n\ge 1$ is fixed, then $G_{2n}=\es$ iff $\g_{n}=\es$. Moreover, 
$|G_{2n}|=o(1)$ as $n\to \iy$.

\no (v) If ${N\/3}\in \Z$, then $
G_{2n+1}=\es$ for all $n\ge 0$ iff $q\in L_{even}(0,1)$. Moreover, 
 each odd gap $G_n$ has the form $G_n=(r_{p,n}^-,r_{p,n}^+)$ ($n\ge 1$ is odd) for some integer $p=p(n)$ and $r_{p,n}^\pm$ satisfy  
\[
\lb{Tas-6}
r_{p,n}^\pm=\pi^2{n^2\/4}+q_0\pm |\wt q_{cn}|+o(n^{-1})\ \  as\qq n\to \iy,
\qq \wt q_{cn}=\int_0^1q(t)\cos \pi nt dt.
\]
\no (vi) If ${N\/3}\notin \Z$, then each gap
$G_{2n+1}\ne \es, n\ge 0$ and $|G_{2n+1}|\to \iy$ as $n\to \iy$.
\end{theorem}

Note that in the armchair case (see \cite{BBKL}) 
there exist non-real resonances for some specific potentials.
That all resonances are  real is a peculiarity of the high
symmetry of a zigzag periodic graph. 

Finally, we formulate our results about {\bf the inverse problem}.
By Theorem \ref{T3}, $E_n=\m_n,n\ge 1$ are the eigenvalues of $H$ and they satisfy (see [PT]) $E_n=\pi^2n^2+q_0+\vk_n$, where $q_0=\int_0^1q(t)dt$ and $(\vk_n)_{1}^\iy\in\ell^2=\ell _0^2$.
Here $\pi^2n^2, n\ge 1$ are the unperturbed eigenvalues and  the real Hilbert spaces $\ell _p^2$ is given by
$
\ell _p^2=\rt\{ f=(f_n)_1^\iy:\sum_{n\ge
1}n^{2p}f_n^2<\iy \rt\}$, $p\ge 0.
$
The monotonicity property $E_1<E_2<...$ gives that if
$q$ runs through $L^2(0,1)$, then  $(\vk_n)_1^\iy$ doesn't run through the whole space $\ell^2$. In order to describe this situation, we introduce the open and convex set
$
\cK=\left\{(\vk_n)_{1}^\iy\in\ell^2: \pi ^2\!+\!\vk_1\!<\! (2\pi )^2\!+\!\vk_2\!<\!\dots\right\}
\ss\ell^2
$.
Let $Hf_n=E_nf_n$ for some nonzero eigenfunction $f_{n}\in \cH(E_n)$. Then $f_{n,\o}\neq 0$ for some $\o\in \cZ$ and
using Theorem \ref{T2} we obtain 
$$
f_{n,\o}(0)=f_{n,\o}(1)=0,\ \ f_{n,\o}'(0)=C_{n,\o},\ \ f_{n,\o}'(1)^2=C_{n,\o}^2e^{2h_n}>0,\ \ h_n\in \R,\ \ any \ \ 
n\ge 1,
$$
for some constant $C_{n,\o}\ne 0$.
 From Theorem \ref{T3}-\ref{T4} we have a simple corollary.

\begin{corollary}\label{T6} 
(i) The operator $H$ has only a finite number of non degenerate gaps $G_n$ iff ${N\/3}\in \Z$ and $q$ is an even
finite gap potential for the operator $-y''+qy$ on the real line.

\no (ii) The mapping
$
\F:q\to \lt(q_0,(\vk_n)_{1}^\iy;(h_n)_{1}^\iy\rt)
$
is a real-analytic isomorphism between $L^2(0,1)$ and $\R\ts\cK\ts\ell^2_1$, where the set
$
\cK=\left\{(\vk_n)_{1}^\iy\in\ell^2: \pi ^2\!+\!\vk_1\!<\! (2\pi )^2\!+\!\vk_2\!<\!\dots\right\}
\ss\ell^2
$.

\no (iii) The mapping
$
\F_e:q\to \lt(q_0,(\vk_n)_{1}^\iy\rt)
$
is a real-analytic isomorphism between $L_{even}^2(0,1)$ and $\R\ts\cK$.

\end{corollary}

\begin{figure}
\noindent
\centering
\tiny
\psfrag{A}[r][r]{$A_1$}
\psfrag{B}[r][r]{$A_2$}
\psfrag{a}{$\O_1$}
\psfrag{b}{$\O_2$}
\psfrag{c}{$\O$}
\psfrag{(1,0)}{$(1,0)$}
\psfrag{(2,0)}{$(2,0)$}
\psfrag{(3,0)}{$(3,0)$}
\psfrag{(4,0)}{$(4,0)$}
\psfrag{(1,1)}{$(1,1)$}
\psfrag{(2,1)}{$(2,1)$}
\psfrag{(3,1)}{$(3,1)$}
\psfrag{(2,2)}{$(2,2)$}
\psfrag{(3,2)}{$(3,2)$}
\includegraphics[width=.5\textwidth]{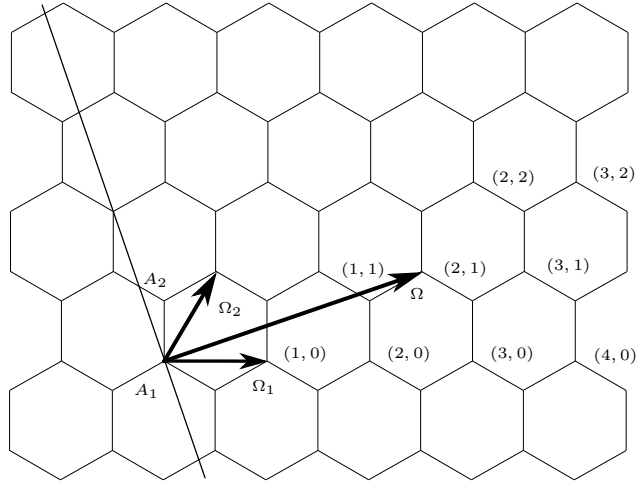}
\caption{The unrolled honeycomb lattice of nanotube.
The unit cell is spanned by the vectors $\O_1$ and $\O_2$. 
The type of the nanotube is defined by the pair $(N_1,N_2)\in\N^2$,
$N_1\geq N_2$, and corresponding chiral vector
$\O=N_1\O_1+N_2\O_2$.
}
\label{fig5}
\end{figure}

For  the convenience of the reader we  briefly describe the structure of carbon nanotubes (see \cite{Ha}, \cite{SDD}) and explain how they are related to the graph $\G^N$. Graphene is a single 2D layer of graphite forming a honeycomb lattice, as in Fig. \ref{fig1} and \ref{fig2}.
A carbon nanotube is a honeycomb lattice "rolled up" into a cylinder.
In carbon nanotubes, the graphene sheet is "rolled up" in such a way 
that a so-called chiral  vector $\O=N_1\O_1+N_2\O_2$ becomes
the circumference of the tube, where $\O_1=e_1+e_2, \O_2=e_0+e_1$
see Fig \ref{fig5} and \ref{fig2}. The chiral  vector $\O$, which is usually denoted by the pair of integers $(N_1,N_2)$, uniquely defines a particular tube. Tubes of type $(N,0)$ are called zigzag tubes, see Fig \ref{fig5} and \ref{fig4}. They exhibit a zigzag pattern along the circumference, see Figure \ref{fig1}. $(N,N)$-tubes are called armchair tubes. 

Remark that in this zigzag 
 standard physical model only each vertical edge has length $1$ and the other edges are shorter (or might be shorter). In this paper, we avoid this additional technical difficulty. Instead, we do  
our detailed  spectral analysis for the simplest possible model
of a physical nanotube $T^N$ (to be described below and see Fig. \ref{fig4}). This will simplify our analysis (in another paper \cite{KL}) of the Schr\"odinger    operator  $\mH_B =(-i\nabla-\mA)^2+q$ on $T^N$ 
with a periodic potential $q$ and  a uniform magnetic field $\mB=B(0,0,1)\in \R^3$, $B\in\R$. The corresponding vector potential is given by $\mA(x)={1\/2}[\mB,x]={B\/2}(-x_2,x_1,0),\  x=(x_1,x_2,x_3)\in\R^3$. Maybe the second simplest model is the armchair nanotube [Ha], but then the monodromy matrix is $4\ts 4$ 
and there exist complex resonances for some specific potential
$q$ (see [BBKL]).

Finally we consider the Schr\"odinger operator $\mH =-\pa^2 +q$
with a periodic potential $q$ on the zigzag nanotube $T^N\ss\R^3$.
Our model nanotube $T^N$ is a union of edges $T_\o$ of length 1, i.e.,
$T^N=\cup_{\o\in \cZ} T_\o$, see Fig. \ref{fig4}. Each edge $T_\o=\{x=
{\bf r}_\o+t{\bf e}_\o,\  t\in [0,1]\}$ is oriented by the vector ${\bf  e}_\o\in \R^3$ and has starting point $\bf r_\o\in \R^3$.
We have the coordinate $x={\bf r}_\o+t{\bf e}_\o$ and the local coordinate $t\in [0,1]$ (length preserving).
We define ${\bf r}_{\o}, {\bf e}_{\o},\o=(n,j,k)\in \cZ$ by
$$
{\bf r}_{n,0,k}=\vk_{n+2k}+{3n\/2}{\bf e}_0,\qq {\bf r}_{n,1,k}={\bf r}_{n,0,k}+{\bf e}_0,\qq {\bf r}_{n,2,k}={\bf r}_{n+1,0,k},
$$
$$
{\bf e}_{n,0,k}={\bf e}_0,\qq {\bf e}_{n,1,k}=
\vk_{n+2k+1}-\vk_{n+2k}+{{\bf e}_0\/2},\qq 
{\bf e}_{n,2,k}=
\vk_{n+2k+2}-\vk_{n+2k+1}-{{\bf e}_0\/2},
$$
where $\vk_j=R_N(c_j,s_j,0)\in \R^3, \ R_N={\sqrt 3\/4\sin {\pi\/2N}}$, 
The points ${\bf r}_{0,0,k}$ are vertices of the regular N-gon $\mP_0$.
The vertical edge $T_{0,0,k}$ is the segment  and  lies on the cylinder $\cC\ev\{(x_1,x_2,x_3)\in \R^3:x_1^2+x_2^2=R_N^2\}
$. The starting points
$
{\bf r}_{1,0,k}={\bf r}_{0,2,k}=\vk_{1+2k}+{3\/2}{\bf e}_0,\ k\in \Z_N
$ 
 are the vertices of the regular N-gon $\mP_1$.
$\mP_1$ arises from $\mP_0$ by the following motion:
rotate around the axis of the cylinder $\cC$ by the angle ${\pi\/N}$
and translate by ${3\/2}e_0$. The non-vertical edges $T_{0,1,k}$ and $T_{0,2,k}$
have  positive and negative projections on the vector ${\bf e}_0$.
Repeating this procedure we obtain all edges of $T^N$. 
Note that each non-vertical edge $T_{0,j,k}, j=1,2$ (without the endpoints) lies inside the cylinder $\cC$.

Our operator $\mH$ on $T^N$ acts
in the space $L^2(T^N)= \os_{\o\in\cZ} L^2(T_\o)$. 
Then acting on the edge $T_\o$, $\mH$ is the ordinary differential operator given by
$(\mH f)_\o=-f_\o''+qf_\o, \ f_\o, f_\o'' \in L^2(\G_\o)$, $\o\in \cZ$
where  $q\in L^2(0,1)$ and  $f\in \gD(\mH)=W^2(T^N)$
satisfies the Kirchhoff Conditions \er{KirC}. Thus the operator
$\mH$ on $T^N$ is unitarily equivalent to the operator $H$
on $\G^N$.

We remark that such models  were introduce by Pauling \cite{Pa}
in 1936 to simulate aromatic molecules. They were described 
in more detail by  ~Ruedenberg and Scherr \cite{RS} in 1953.
For further physical models see [Ha], [SDD].
For papers on spectral analysis of such operators, see the references in \cite{Ku}. 

We present the plan of the paper. In Sect. 2 we construct 
the fundamental solutions and describe the basic properties of the monodromy operator. Moreover, we prove Theorem 1.2 about the eigenfunctions. In Sect. 3 we prove the basic properties of the 
Lyapunov function $\D_0$ which is entire on the plane.
This function is important to study the functions $\D_k, k\in \Z_N$ which are analytic on the two-sheeted Riemann
surface $\cR_k$. In Sect. 4 we prove the basic properties of the 
Lyapunov function $\D_k, k\in \Z_N$ and Corollary \ref{T6}.

\section{Fundamental solutions and eigenfunctions}
\setcounter{equation}{0}

\no {\bf Proof of Theorem \ref{T1}.} (i)
Define the operator $\cS$ in $\C^N$ by
$\cS u=(u_N,u_1,\dots,u_{N-1})^\top$, $u=(u_n)_1^N\in \C^N$.
The unitary operator $\cS$ has the form
$\cS =\sum_1^Ns^k\cP_k$, where $\cS v_k=s^kv_k$ and $v_k={1\/N^{1\/2}}(1,s^{-k},s^{-2k},...,s^{-kN+k})$
 is an eigenvector (recall $s=e^{i{2\pi \/N}}$); 
 $\cP_ku=v_k(u,v_k)$ is a projector. 
The function $f$ in the Kirchhoff boundary conditions \er{KirC} is a vector function  $f=(f_{\o}), \o=(n,j,k)\in\cZ$.  We define a new vector-valued function $f_{n,j}=(f_{n,j,k})_{k}^{N}$, where each $f_{n,j},(n,j)\in \cZ_1=\Z\ts \{0,1,2\}$ is an $\C^N-$ vector, which satisfies the equation
$-f_{n,j}''+qf_{n,j}=\l f_{n,j}$,
and the conditions (which follow from 
the Kirchhoff conditions \er{KirC})
\[
\label{C0}
f_{n,0}(1)=f_{n,1}(0)=\cS f_{n,2}(1),\quad
f_{n+1,0}(0)=f_{n,1}(1)=f_{n,2}(0),
\]
\[
\label{C1}
-f'_{n,0}(1)+f'_{n,1}(0)-\cS f'_{n,2}(1)=0,\quad
f'_{n+1,0}(0)-f'_{n,1}(1)+f'_{n,2}(0)=0,
\]
for all $n\in\Z$. The operators $\cS$ and $H$ commute, then
$H=\os_1^N (H\cP_k)$.
Using \er{C0}, \er{C1}  we deduce that $H\cP_k$ is unitarily equivalent to the operator $H_k$. The operator $H_k$ on $\G^1=\cup_{\a\in\cZ_1} \G_\a^1$ acts in the Hilbert space $L^2(\G^1)$. 
Then acting on the edge $\G_\a$,  the operator $H_k$ given by:
$ (H_kf)_\a=-f_\a''+q(t)f_\o$, where $f_\a,f_\a''\in L^2(0,1),
\a=(n,j)\in\cZ_1$,
on the vector functions $f=(f_\a)_{\a\in\cZ_1}$, which 
satisfy the boundary conditions \er{1K0}, \er{1K1}, i.e.,
\[
\label{K0}
f_{n,0}(1)=f_{n,1}(0)=s^k f_{n,2}(1),\qq
f_{n+1,0}(0)=f_{n,1}(1)=f_{n,2}(0), \qq s=e^{i{2\pi\/N}},
\]
\[
\label{K1}
-f'_{n,0}(1)+f'_{n,1}(0)-s^k f'_{n,2}(1)=0,\qqq
f'_{n+1,0}(0)-f'_{n,1}(1)+f'_{n,2}(0)=0.
\]

(ii) Fix $k\in \Z_N$.
 We determine the fundamental solutions
$\vT_k(x,\l),\F_k(x,\l)$ of the equation $-f''+qf=\l f,\qq \l\in\C$,
on $\G^1$, where $\vT_k,\F_k$  satisfy \er{K0}-\er{K1} and  
\[
\lb{inco}
\vT_{k,\a}(0,\l)=\F_{k,\a}'(0,\l)=1,\qq 
\vT_{k,\a}'(0,\l)=\F_{k,\a}(0,\l)=0,\qqq \a=(0,0).
\]
Below we assume $f=\vT_k$ or $f=\F_k$ and let $\vt_t=\vt(1,\l), \vp_t=\vp(t,\l),...$.
Any solution of the equation $-f''+qf=\l f$ satisfies
$f(t)=\vt_tf(0)+{\vp_t\/\vp_1}(f(1)-\vt_1f(0)),\qq t\in[0,1]$
and 
\[
\label{gs2}
f'(0)={f(1)-\vt_1f(0)\/\vp_1},\qq
f'(1)={\vp_1'f(1)-f(0)\/\vp_1}.
\]
Substituting \er{gs2} for $f=f_{01}, f=f_{02}$ into the first Eq. in  \er{K1} at $n=0$, we obtain
$$
-\vp_1f'_{0,0}(1)+(f_{0,1}(1)-\vt_1f_{0,1}(0))
-s^k(\vp_1'f_{0,2}(1)-f_{0,2}(0))=0.
$$
Then the condition \er{K0} gives
\[
\vp_1f'_{0,0}(1)+( f_{1,0}(0)-\vt_1 f_{0,0}(1))
-s^k(\vp_1's^{-k}f_{0,0}(1)-f_{1,0}(0))=0,
\]
which implies 
\[\label{sing}
-\vp_1f'_{0,0}(1)+u f_{1,0}(0)-2\D f_{0,0}(1)=0,\qqq
u=1+s^k=2c_ks^{{k\/2}}.
\]
Thus we obtain \textbf{the first basic identity}
\[
\lb{t1}
f_{1,0}(0)=
{\vp_1\/u} f'_{0,0}(1)+{2\D\/u} f_{0,0}(1),\qqq u=2c_ks^{{k\/2}},\qq
c_k=\cos {\pi k\/ N}.
\]
Substituting \er{gs2} for $f=f_{01}, f=f_{02}$ into the second Eq. in \er{K1} at $n=0$, we get
$$
\vp_1f'_{1,0}(0)-(\vp_1'f_{0,1}(1)-f_{0,1}(0))
+(f_{0,2}(1)-\vt_1f_{0,2}(0))=0.
$$
Using \er{K0}, we obtain
\[
\vp_1f'_{1,0}(0)-(\vp_1'f_{1,0}(0)-f_{0,0}(1))
+(s^{-k}f_{0,0}(1)-\vt_1f_{1,0}(0))=0,
\]
which yields
\[\label{t2}
\vp_1f'_{1,0}(0)=2\D f_{1,0}(0)-u_2 f_{0,0}(1),\qqq u_2=1+s^{-k}=2c_ks^{-{k\/2}}.
\]
Then substituting \er{t1} into \er{t2}, we obtain 
\textbf{the second basic identity}
\[
\label{psi31}
f_{1,0}'(0)={4(\D^2-c_k^2)\/u\vp_1}f_{0,0}(1)+
{2\D\/u} f_{0,0}'(1), \qqq u=2c_ks^{{k\/2}}.
\]
Substituting consequently $f=\vT_k$ and $f=\F_k$ into
\er{t1} and into \er{psi31}, and after this $\vT_{k,\o}(0,\l), \vT_{k,\o}'(0,\l),..$ into 
$\cM_k(\l)=\ma \vT_{k,\a}(0,\l) & \F_{k,\a}(0,\l)\\
\vT_{k,\a}'(0,\l) & \F_{k,\a}'(0,\l)\am,\ \a=(1,0)$, we obtain
\er{T1-1}.

Using \er{T1-1} we obtain 
$
\det \cT_k={s^{-k}4c_k^2\/4c_k^2} =s^{-k}, \  \det \cM_k=\det \cT_k\cM=s^{-k}$, 
since $\det \cM=1$. Moreover, using  
$
\cR \cM\cR^{-1}=\ma \vt_1&1\\\vp_1\vt_1'&\vp_1'\am, \qq
\ \ \cR=\ma 1 & 0\\0 & \vp_1\am,
$
we obtain 
$$
\Tr\cM_k=\Tr\cT_k \cR \cM\cR^{-1}={2s^{-{k\/2}}\/c_k}\Tr \ma 2\D & 1 \\ 4\D^2-4c_k^2 & 2\D \am
\ma \vt_1&1\\
\vp_1\vt_1'&\vp_1'\am\\
$$$$
={2s^{-{k\/2}}\/c_k}\rt(2\D \vt_1+\vp_1\vt_1'+4(\D^2-c_k^2)+2\D \vp_1'\rt)
={2s^{-{k\/2}}\/c_k}
\rt(2\D^2+{\vp_1\vt_1'\/4}-c_k^2\rt)={2s^{-{k\/2}}(\D_0+s_k^2)\/c_k}
$$
which gives \er{T1-2}.
\BBox

Recall that $H^0, \gD(H^0)=W^2(\G^N)$ is  self-adjoint [Ca1].
Let $\|f\|^2=\int_{\G^N}|f(x)|^2dx$ and let $\|f_\o\|_0^2=
\int_0^1|f_\o(t)|^2dt$ for $f_\o\in L^2(\G_\o),\o\in \cZ$.
 Assume that 
\[
\lb{selfa}
\|qf\|^2\le {1\/2}\|H^0f\|^2+C\|f\|^2, \ \ all\ \ f\in \gD(H^0),
\]
for some constant $C>0$. Then by the Kato-Rellich Theorem (162 p. [ReS1]), $H=H^0+q$ is self-adjoint on $\gD(H^0)$
and essentially self-adjoint on  any core of $H^0$.
We prove \er{selfa}.

\begin{lemma}
\label{self-est} 
Let $q\in L^2(0,1)$. Then the estimate \er{selfa} holds true.
\end{lemma}

\no {\it Proof.} For $f\in \gD(H^0), f=(f_\o)_{\o\in \cZ}, \o=(n,j,k)\in \cZ$, we have the following estimates 
$$
\|qf_\o\|_0 \le \max_{x\in \G_\o} |f_\o(x)|\|q\|_0,\qq
\max_{x\in \G_\o} |f_\o(x)|\le \int_{\G_\o}(|f_\o|+|f_\o'|)dx\le 
\|f_\o'\|_0+\|f_\o\|_0,
$$
which yield
$$
\|qf\|^2=\sum \|qf_\o\|_0^2\le \|q\|_0^2\sum (\|f_\o'\|_0+\|f_\o\|_0)^2.
$$
Thus for $C_q=\|q\|_0^2$ we obtain
$$
\|qf\|^2\le 2C_q\sum (\|f_\o'\|_0^2+ \|f_\o\|_0^2)=2C_q(\|f'\|^2+
\|f\|^2)\le
\ve\|H^0f\|^2+\rt({C_q^2\/\ve}+2C_q\rt)\|f\|^2 ,\ \ 
$$
for any $\ve>0$, since 
$$
2C_q\|f'\|^2=2C_q(H^0f,f)\le 2C_q\|H^0f\|\|f\|\le
\ve\|H^0f\|^2+{C_q^2\/\ve}\|f\|^2,
$$
which yields  \er{self-est}. \BBox

Note that the operator $H_k$ in $L^2(\G^1)$
with $\gD(H_k)=\{f, \os_{\a\in \cZ_1}f_\a''\in L^2(\G^1): f\ \   
satisfies $ $\ the \ Kirchhoff\  conditions \ \er{1K0}, \ \er{1K1} \}$
is self-adjoint. The proof is similar to the case $H$.

\no {\bf Proof of Theorem \ref{T2}}.
(i) Let $H_kF=\l F$ for some eigenfunction $F$ and some $\l\in\s_D$.
For each  $F_\o=F|_{\G_\o}, \o=(n,j)\in\cZ_1=\Z\ts \J$,
we have $F_\o(t)=F_\o(0)+\int_0^t F_\o'(s)ds$, $t\in[0,1]$. Then
$|F_\o(0)|\leq |F_\o(x)|+\int_0^1|F_\o'(t)|dt
$ and integrating over $x\in[0,1]$, using the H\"older inequality, we obtain 
\[\label{Fto0}
|F_\o(0)|\le \int_0^1 (|F_\o(t)|+|F_\o'(t)|)dt\le \|F_\o\|_0+\|F_\o'\|_0=o(1), \qq as \qq n\to\pm\iy,
\]
since $F\in\cH_k(\l)$. Furthermore, each restriction $F_\o$ has 
the form $F_\o(t)=a_\o\vt_t+b_\o\vp_t, t\in [0,1]$
for some constants $a_\o, b_\o$, which implies  $F_\o(1)=F_\o(0)\vt_1$, where $\vp_t=\vp(t,\l),\vt_t=\vt(t,\l),..$
The Kirchhoff conditions give
$F_{n+1,0}(0)=F_{n,1}(1)=\vt_1F_{n,1}(0)=\vt_1F_{n,0}(1)=\vt_1^2F_{n,0}(0)$, which yields 
\[
\label{Fninf}
F_{0,0}=\vt_1^{2n}F_{-n,0}(0),\ \ \ 
all \ n\in\Z.
\] 
Let $|\vt_1|\le 1$, the proof for the case
$|\vt_1|> 1$ is similar. Then \er{Fninf}, \er{Fto0} imply $F_{0,0}=\vt_1^{2n}F_{-n,0}(0)=o(1)$ as  $n\to+\iy$.
Thus \er{Fninf} gives $F_{n,0}(1)=0$ for all $n\in\Z$.
Finally, $F$ vanishes on all vertexes of $\G^1$, since the set of all ends of vertexes $\G_{n,0}, n\in\Z$ coincides with the vertex set of $\G^1$.

(ii) 
Using \er{T2-1}-\er{T2-2}, we deduce that  $\p^{(0)}$ satisfies
the Kirchhoff conditions \er{1K0}, \er{1K1}. Thus
$\p^0$ is a eigenfunction  of $H_k$.

The operator $H_k$ is periodic, then each $\p^{(n)}, n\in \Z$ is an eigenfunction. We will show that the sequence $\p^{(n)}, n\in\Z$ forms
a basis for $\cH_k(\l)$. 
The functions $\p^{(n)}$ are linearly independent, since $\supp\p^{(n)}\cap\supp\p^{(m)}=\es$  for all $n\ne m$.

{\bf Consider the first case} $\e=1-s^k{\vp_1'}^2\ne 0$. 
For any $f\in \cH_k(\l)$ we will show the identity \er{T2-3}, i.e.,
$f=\wh f$, where $\wh f=\sum_{n\in \Z}\wh f_n\p^{(n)},\ \wh f_n={f_{n,0}'(0)\/\e}$. The definition of $\wh f$ and $\l\in \s_D$ give
$\wh f|_{\G_{n,0}}=f|_{\G_{n,0}}=\wh f_n\vp$ for all $n\in\Z$.
This yields $\sum |\wh f_n|^2<\iy$ and $\wh f\in L^2(\G^1)$, since $f\in L^2(\G^1)$

Note that $\wh f$ satisfies the Kirchhoff conditions \er{1K0}, \er{1K1}
and $-\wh f_\a''+q\wh f_\a=\l \wh f_\a, \a\in \cZ_1$.
Consider the function $u=f-\wh f=(u_\a)_{\a\in \cZ_1}$, where
$u_{n,j}=C_{n,j}\vp(t,\l),  j=1,2$ and $u_{n,0}=0$ for some constant $C_{n,j}, n\in \Z$. The Kirchhoff boundary conditions \er{1K0}-\er{1K1} yields $u=0$.

{\bf Consider the second case} $s^k{\vp_1'}^2=1$.
For any $f\in \cH_k(\l)$ we will show the identity \er{T2-3}, i.e.,
$f=\wh f$, where
$
\qq \wh f=\sum_{n\in \Z}\wh f_n\p^{(n)},\ \wh f_n=f_{n,1}'(0).
$
From the definition of $\wh f$ and $\l\in \s_D$ we deduce that
$\wh f|_{\G_{n,1}}=f|_{\G_{n,1}}=\wh f_n\vp $ for all $n\in\Z$.
This yields $\sum |\wh f_n|^2<\iy$ and $\wh f\in L^2(\G^1)$, since $f\in L^2(\G^1)$.

Consider the function $u=f-\wh f$. 
The function $u=0$ at all vertices of $\G^{1}$
and then $u_\a=C_\a\vp_t, \a=(n,j), n\in \Z, j=0,2$.
Assume that $C_{0,0}=C$. Then the Kirchhoff boundary conditions
\er{1K1} yields $C_{n,0}=-C_{n,2}=-C$ and $C_{n,0}=C_{n+1,0}=C$,
which give $C=0$ since $u\in L^2(\G^{(1)})$.
\BBox

\section{ The Lyapunov function $\D_0$}
\setcounter{equation}{0}

Recall that $\vt(x,\l), \vp(x,\l), x\in \R$ are the fundamental solutions  of Eq. $-y''+qy=\l y,\ \l\in\C,$ on the real line
such that $\vt(0,\l)=\vp'(0,\l)=1,\  \vt'(0,\l)=\vp(0,\l)=1.
$ For each $x\in \R$ the functions $\vt,\vt',\vp,\vp'$ are entire in $\l\in\C$ and satisfy:
\[
\lb{vtas} 
\vt(x,\l)=\cos\sqrt{\l}x+
{1\/2\sqrt{\l}}\int_0^x\lt(\sin\sqrt{\l}x+\sin\sqrt{\l}(x-2t)\rt)q(t)dt+O\lt({e^{|\Im\sqrt{\l}|x}\/|\l|}\rt),
\]
\[
\lb{vt'as} \vt'(x,\l)=-\sqrt{\l}\sin\sqrt{\l}x+
{1\/2}\int_0^x\left(\cos\sqrt{\l}x+\cos\sqrt{\l}(x-2t)\right)q(t)dt+
O\lt({e^{|\Im\sqrt{\l}|x}\/|\l|^{1/2}}\rt),
\]
\[
\lb{vpas} \vp(x,\l)= {\sin\sqrt{\l}x\/\sqrt{\l}}+
{1\/2\l}\int_0^x\left(-\cos\sqrt{\l}x+\cos\sqrt{\l}(x-2t)\right)q(t)dt+
O\lt({e^{|\Im\sqrt{\l}|x}\/|\l|^{3/2}}\rt),
\]
\[
\lb{vp'as} \vp'(x,\l)= \cos\sqrt{\l}x+
{1\/2\sqrt{\l}}\int_0^x\left(\sin\sqrt{\l}x-\sin\sqrt{\l}(x-2t)\right)q(t)dt+
O\lt({e^{|\Im\sqrt{\l}|x}\/|\l|}\rt)
\]
as $|\l |\to\iy$, uniformly on bounded sets of $(x;q)\in [0,1]\ts L^2(0,1)$  (see [PT]), and
\[
\label{Das} \D(\l)= \cos\sqrt{\l}+
{q_0\sin \sqrt{\l}\/2\sqrt{\l}}+{O(e^{|\Im\sqrt{\l}|})\/|\l|},\ \ 
q_0=\int_0^1q(t)dt.
\]
Substituting \er{vtas}-\er{Das} into $\D_0$ we obtain \er{asD0},
i.e., 
$$
 \D_0(\l)=\D_0^0(\l)+{9q_0\/8}{\sin 2\sqrt{\l}\/\sqrt{\l}}+O\lt({e^{2|\Im\sqrt{\l}|}\/|\l|}\rt),\ \ 
\D_0^0(\l)={9\cos 2\sqrt{\l}-1\/8}.
$$
Let $\D_-={1\/2}(\vp'(1,\cdot)-\vt(1,\cdot))$. Substituting the identity
$
\vp(1,\cdot)\vt'(1,\cdot)=\D^2-\D_-^2-1
$
into $\D_0=2\D^2(\l)+{\vt'(1,\l)\vp(1,\l)\/4}-1$, we obtain
\[
\lb{DY}
\D_0={9\D^2-\D_-^2-5\/4},\ \qq {\rm where} \qq \D_-={\vp'(1,\cdot)-\vt(1,\cdot)\/2}.
\]
Asymptotics \er{vtas}-\er{vp'as} yield as $ |\l|\to \iy$
\[
\lb{asY}
\D_-(\l)=-{F(\l)\/2\sqrt \l}+{O(e^{|\Im\sqrt{\l}|})\/|\l|},\ \ \qq
F(\l)=\int_0^1\sin \sqrt \l(1-2t)q(t)dt.
\]

\no {\bf Proof of Theorem \ref{T3}.}
(i) Using \er{T1-2} we obtain the characteristic equation 
$$
\det(\cM_k-\t I_2)=\t^2-2a_k\t+s^{-k}=0, \ \ \ \ a_k={\Tr \cM_k\/2}
={2\x_k\/1+s^k},\qq \x_k=\D_0+s_k^2.
$$
The eigenvalues $\t_{k,\pm}$ of $\cM_k$ are given by
$
\t_{k,\pm}=a_k\pm\sqrt{a_k^2-s^{-k}}.
$
Using $a_k{(1+s^{k})\/2}=\x_k$ and ${(1-s^{k})^2\/4}(a_k^2-s^{-k})=\r_k$  we have
$$
\D_{k,+}={1\/2}(\t_{k,+}+\t_{k,+}^{-1})=
{1\/2}(\t_{k,+}+s^{k}\t_{k,-})=a_k{(1+s^{k})\/2}+
{(1-s^{k})\/2}\sqrt{a_k^2-s^{-k}}=\x_k+\sqrt{\r_k}
$$
and a similar argument yields $\D_{k,-}=\x_k-\sqrt{\r_k}$.

(ii) The standard arguments (see [Ca1]) yields
$\s(H_k)=\s_{ac}(H_k)\cup\s_{pp}(H_k)$, i.e., $\s_{sc}(H_k)=\es$
and $\s_{ac}(H_k)=\{\l\in\R: \D_k(\l)\in [-1,1]\}$.

(iii) Let $\D_k'(\l_0)=0$
and $\D_k(\l_0)\in (-1,1)$ for some $\l_0\in \s_k, k=0,..,N$. Then we have the Tailor series
$
\D_k(\l)=\D_k(\l_0)+\z^p{\D_k^{(p)}(\l_0)\/p!}+O(\z^{p+1})$ as $\z=\l-\l_0\to 0$,
where $\D_k^{(p)}(\l_0)\ne 0$ for some $p\ge 1$. By the
Implicit Function Theorem, there exists some curve $Y\ss
\{\l:|\l-\l_0|<\ve \}\cap \C_+, Y\neq \es$,  for some $\ve>0$
such that $\D_k(\l)\in (-1,1)$ for any $\l\in Y$. Thus we have
a contradiction with \er{T3-1}.
\BBox

In order to prove Theorem \ref{T4} we need

\begin{lemma}
\label{T41} Let $q\in L^2(0,1)$. Then the following statements hold:

\no (i) $\D_0(\l)\ge 1$ for all $\l=\m_n$ and $\l=\n_{n-1},n\ge 1$.
If in addition $\D^2(\l)=1$, then $\D_0(\l)=1$.

\no (ii) If  $\D(\e)=0$ for some $\e\in \R$, then $\D_0(\e)=-{5+\vt^2(1,\e)\/4}\le -{5\/4}$. 

\no (iii) Let $q\in L^2_{even}(0,1)=\{q\in L^2_{even}(0,1):q(x)\!=\!q(1\!-\!x), x\!\in\![0,1]\}$. Then
\[
\D_0={9\D^2-5\/4},\ \ \ \ \D=\vp'(1,\cdot)=\vt(1,\cdot),\ \  
\]
\[
\g_{0,2n}=\g_n,\ \ where\ each \ 
\g_{n}=(\m_n,\n_n),\ \  or\ \ \ \g_{n}=(\n_n,\m_n),
\ \ \ \ n\ge 1.
\]
The zeros $\e_n=\l_{0,2n-1}$ for all $n\ge 1$.

\no (iv) There exists an integer $n_0>1$ such that 
 $\D_0-1$ has exactly $2n_0+1$ roots,  counted with multiplicities, in the domain $\{\l: |\sqrt{\l}|<\pi (n_0+{1\/2})\}$
and for each $n>n_0,$ exactly two roots,  counted with multiplicities, in the domain
$\{\l: |\sqrt \l-\pi n|<{\pi\/4}\}.$  There are no other roots.

\no (v) There exists an integer $n_0>1$ such that 
 $\D_0+{5\/4}$ has exactly $2n_0$ roots,  counted with multiplicities, in the domain $\{z: |\sqrt{z}|<\pi n_0\}$
and for each $n>n_0,$ exactly one simple root in the domain
$\{\l: |\sqrt \l-\pi (n-{1\/2})|<{\pi\/4}\}.$  There are no other roots.

\no (vi) There exists an integer $n_0>1$ such that 
 $\D_0$ has exactly $n_0$ roots,  counted with multiplicities, in the domain $\{\l: |\sqrt{\l}|<\pi n_0\}$
and for each $n>n_0,$ exactly one simple root in the domain
$\{\l: |\sqrt \l-\pi (n-{1\/2})|<{\pi\/4}\}.$  There are no other roots.

\no (vii) Let $c\in (-{5\/4},1)$
and let $u_n^{\pm}=\pi n\pm u_0^{+}, u_0^{+}={\arccos {1+8c\/9}\/2}\in [0,{\pi\/2}],\ n\ge 1$ and $u_{0}^{+}<u_{1}^{-}
<u_{1}^{+}<u_{2}^{-}<.. $. Then there exists an integer $n_0>1$ such that  $\D_0-c$ has exactly $2n_0+1$ roots,  counted with multiplicities, in the domain $\{\l: |\sqrt{\l}|<R\},
R=(u_{n_0}^{+}+u_{n_0+1}^{-})/2$
and for each $n>n_0,$ exactly one simple root in the domain
$\{\l: |\sqrt \l-u_{n}^\pm|<r$ for some small $r\in (0,1)$.  There are no other roots.
\end{lemma}

\no{\it Proof.} (i) Let $\l\in\{\m_n,\n_n\}$ for some $n\ge 1$.
Then we have $\D_0(\l)=2\D^2(\l)+{\vt'(1,\l)\vp(1,\l)\/4}-1=2\D^2(\l)-1  \ge 1$. Moreover, if $\l\in\{\m_n,\n_n\}$ and $\D^2(\l)=1$, then we have $\D_0(\l)=1$.

 (ii) We have $0=2\D(\e)=\vp'(1,\e)+\vt(1,\e)$,
which yields $\vp'(1,\e)=-\vt(1,\e)$. Substituting the last
identity into $\D_0(\l)=2\D^2(\l)+{\vt'\vp\/4}-1$ we get the
needed estimate.

(iii)  It is well known that if $q\in L^2_{even}(0,1)$,
then $\D=\vp'(1,\cdot)=\vt(1,\cdot)$, see p.8, [MW].
It is well known that if $\D=\vp'(1,\cdot)=\vt(1,\cdot)$, then $q\in L^2_{even}(0,1)$, see [PT].

(iv)  Define the contours $C_n(r)=\{\l\in \C:|\sqrt \l-\pi n|=\pi r\}$. Fix another integer $n_1>n_0$. Consider the contours
$C_{0}(n_0+(1/2)), \ C_{0}(n_1+(1/2))),\ C_{n}(1/4), n>n_0.$
Using the estimate $e^{|\Im \sqrt\l|}<4|\sin \sqrt\l|,|\l-\pi n|\ge {\pi \/4}$ and \er{asD0} we obtain on all contours (for large $n_0$)
$$
|(\D_0(\l)-1)-(\D_0^0(\l)-1)|=o(e^{2|\Im \sqrt \l|})=
|\sin \sqrt \l|^2o(1)=o(1)|\D_0^0(\l)-1|
$$
Hence, by Rouch\'e's Theorem, $\D_0(\l)-1$ has as many roots, counted with multiplicities, as $\D_0^0(\l)-1$ in each of the bounded domains and the remaining unbounded domain. Since $\D_0^0(\l)-1$ has exactly the simple root $=0$ and one root 
of the multiplicity $2$ at each $\pi^2n^2\ge 1$, and since $n_1>n_0$ can be chosen arbitrarily large, the point (iii) follows.

 The proof of (v) and (vi) is similar.
 \BBox

\begin{lemma}
\label{TLas} Let $c\in [-{5\/4},1]$. Then the equation $\D_0(\l)=c,\l\in \C$
has only real zeros, which satisfy

\no (i) If $c\in (-{5\/4},1)$ and $q_0=\int_0^1q(t)dt$, then these zeros $z_n^\pm$ are given by
\[
\lb{TLas-1}
z_{0}^{+}<z_{1}^{-}<z_{1}^{+}<z_{2}^{-}<.. \qq and \qq
\sqrt{z_{n}^{\pm}}=u_{n}^{\pm}+{q_0\/2u_n^{\pm}}+{o(1)\/n^2}\ \ as \ \ \ \ n\to \iy,
\]
where $u_n^{\pm}=\pi n\pm u_0^{+}, \ n\ge 1$ and $u_{0}^{+}={\arccos {1+8c\/9}\/2}\in [0,{\pi\/2}],u_{0}^{+}<u_{1}^{-}
<u_{1}^{+}<u_{2}^{-}<.. $.

Moreover, the zeros $z_n^\pm$ have  another forms given by: $x_n^-=z_{n+1}^{+}, x_n^+=z_{n}^{-},$  and
\[
\lb{TLas-2}
x_{1}^{-}<x_{1}^{+}<x_{2}^{-}<x_{2}^{+}<x_{3}^{-}<.. \qqq
\sqrt{x_{n}^{\pm}}=v_{n}^{\pm}+{q_0\/2v_n^{\pm}}+{o(1)\/n^2}\qq as \qq n\to \iy, 
\]
where $v_n^{\pm}=\pi (n-{1\/2})\pm ({\pi\/2}-u_0^{+})$.

\no (ii) If $c=1$ or $c=-{5\/4}$ . Then these zeros $z_n^\pm$ and $x_n^\pm$ are given by
$$
z_{0}^{+}<z_{1}^{-}\le z_{1}^{+}<z_2^-\le z_2^+<.. \qq and \qq
\sqrt{z_{n}^{\pm}}=\pi n+{q_0\/2\pi n}+{o(1)\/n^2},\qq if \qq
c=1,
$$
$$
x_{1}^{-}\le x_{1}^{+}<x_2^-\le x_2^+<.. \qq and \qq
\sqrt{x_{n}^{\pm}}=\pi (n-{1\/2})+{q_0\/2\pi n}+{o(1)\/n^2},\qq if \qq c=-{5\/4},
$$
as $n\to \iy$.
\end{lemma}

\no {\it Proof.} (i)
Let $z=\sqrt \l$. Recall the folowing asymptotics from [KK1]:
\[
\lb{asLq}
\D(\l)=\cos \c(\l),\qq \c(\l)=\sqrt \l-{q_0\/2\sqrt \l}+{o(1)\/\l} \qqq 
as \qq |\l|\to \iy.
\]
Using \er{DY} we rewrite the Eq. $\D_0(\l)=c,\l\in \R$ in the form
$c={9\D^2-\D_-^2-5\/4}$, which yields 
\[
\lb{asqu}
\cos 2\c(\l)=A_c+{2\D_-^2(\l)\/9},\ \ A_c={1+8c\/9}\in (-1,1).
\]
In the case $q=0$ we have $\cos 2z=A_c$ and the corresponding solutions $\sqrt{z_n^\pm}=u_n^\pm$. Then by Lemma \ref{T41},\
$|\sqrt{z_n^{\pm}}-u_n^{\pm}|<r$ for all $n>n_0,$ for some $n_0\ge 1$
and some small $r\in (0,1)$.  Using \er{asqu},\er{asY} we obtain 
$
\cos 2\c(\l)=A_c+o(\l^{-1})$ as $\l\to \iy$,
and then $\z_n^\pm=\sqrt{z_{n}^{\pm}}-u_{n}^{\pm}\to 0$.
Thus $\ve_n^\pm=\z_n^\pm-{q_0\/2u_{n}^{\pm}}+{o(1)\/n^2}\to 0$
 and the Taylor series gives
$$
\cos 2\c(\l)
=\cos 2(u_n^{\pm}+\ve_n^\pm)=\cos 2u_n^{\pm}-
\sin (2u_n^{\pm})\ve_n^\pm(1+O(\ve_n^\pm))=A_c+o(n^{-2}).
$$
Thus we obtain $\ve_n^\pm=\z_n^\pm-{2_0\/2u_{n}^{\pm}}+{o(1)\/n^2}={o(1)\/z^2}$
which implies, $\z_n^\pm={q_0\/2u_{n}^{\pm}}+{o(1)\/n^2}$.

Proof of \er{TLas-2} and  (ii) is similar.
\BBox

{\bf Proof of Theorem \ref{T4}}.
 (i) The function $\D_0=2\D^2+{\vt'(1,\cdot)\vp(1,\cdot)\/4}-1$ is entire and real on real line. Moreover, it has asymptotics \er{asD0}. By Theorem \ref{T3}, the function $\D_0$ has only simple real zeros $\e_{0,n},n\ge 1$, which satisfy $\e_{0,1}<\e_{0,2}<\e_{0,3}<...$.
Then by the Laguerre Theorem (see Sect. 8.52 [Ti]), the function $\D_0'$ has only real simple zeros $\l_{0,n},n\ge 1$,  which are separated by the zeros of $\D_0$:
$\e_{0,1}<\l_{0,1}<\e_{0,2}<\l_{0,2}<\e_{0,3}<...$,
since $\D_0(\l)\to \iy$ as $\l\to -\iy$. Moreover, we obtain
$(-1)^n\D_0(\l_{0,n})\ge 1,n\ge 1$.

Define an interval $\cG=(-\iy, \pi^2n_0^2)$. Lemma \ref{T41} gives
for large integer $n_0>1$ : 

A) The function $\D_0-1$ has $2n_0+1$ zeros on $\cG$.

B) $\D_0(\e_{0,n})\le -{5\/4}, n=1,..,n_0$, where $\e_{n}\in \cG$ 
is a zero of $\D$, $\e_1<\e_2<...$.

C) $\D_0(\m_n)\ge 1, n=1,..,n_0$, where $\m_n\in \cG$ 
are zeros of $\vp (1,\l)$, the Dirichlet eigenvalues.

These facts and $\m_1<\e_1<\m_2<\e_2<...$  yield 
$\m_n\in \wt\g_{0,2n}$ and $\e_n\in \wt\g_{0,2n-1}$
for all $n\ge 1$, where $\wt\g_{0,n}=[\l_{0,n}^-,\l_{0,n}^+]$. Moreover, we obtain \er{T4-1}.

(ii) Using (i) we deduce that the periodic and anti-periodic eigenvalues $\l_{0,n}^\pm,n\ge 0$ satisfy \er{epa0}.
Consider $\l_{0,2n+1}^\pm$, which  satisfies $\D_0(\l)=-1$.
Using Lemma \ref{TLas} (i) for $c=-1$ we have the first asymptotics in \er{T4-2}. In order to show \er{T4-2} we 
consider $\l_{0,2n}^\pm$ with even $n\ge 1$.
The proof for odd $n$ is similar. 
Using \er{vtas}-\er{vp'as} and \er{asY} we obtain
\[
\lb{asDY}
\!\!\dot \D_0(\l)={O(1)\/n^2},\ \ \ddot \D_0(\l)=-{1+{O(1)\/n}\/(\pi n)^2},\ \
\D_-(\l)={q_{sn}+{O(1)\/n}\/2\pi n},\ \ 
\dot \D_-(\l)=-{\wt q_{cn}+{O(1)\/n}\/(2\pi n)^2},\ \ 
\]
\[
\lb{asD23}
\dot \D(\l)={O(1)\/n^2},\
\ddot \D(\l)=-{1+{O(1)\/n}\/(\pi n)^2},\qq
\dddot \D(\l)={O(1)\/n^4}, \qq where \qq \dot \D={\pa \D\/\pa \l}
\]
as $\sqrt {\l}=\pi n+O(1/n)$,
where $\wt q_{cn}=\int _0^1(1-2t)q(t)\cos 2\pi ntdt$.
Using $\dot \D_0(\l_{0,2n})=0$, we have 
$$
\dot \D_0(\l_{n})=\ddot \D_0(\l_{0,2n})s_{n}(1+O(s_{n})),\qq s_n=\l_{n}-\l_{0,2n} \qq as \qq n\to \iy.
$$
The identity $\dot \D(\l_{n})=0$ and $\D_0={9\D^2-\D_-^2-5\/4}$ 
and \er{asDY} yield
\[
\lb{asD'}
\dot \D_0(\l_{n})=-{\D_-(\l_{n})\dot \D_-(\l_{n})\/2}=
{(q_{sn}+O({1\/n}))(\wt q_{cn}+O({1\/n}))\/(2\pi n)^3}.
\]
Thus asymptotics \er{asDY}-\er{asD'} give
\[
s_n=\l_{n}-\l_{0,2n}=-{(q_{sn}+O({1\/n}))(\wt q_{cn}+O({1\/n}))\/8\pi n},
\]
which implies that roughly speaking the point $\l_{0,2n}$ is in the center of the gap $\g_n$, since $\l_n=(\pi n)^2+q_0+O(1/n)$
see [Ko]. 

We will determine asymptotics of $\g_{0,2n}=(\l_{0,2n}^-,\l_{0,2n}^+)$.
Due to \er{asD23} we get
\[
\D(\l_{0,2n}^\pm)=\D(\l_{n})+A_{0,n}^\pm,\qq  A_{0,n}^\pm=\ddot \D(\l_{n}) {{\z_{0,2n}^\pm}^2\/2}(1+{O(\z_{0,2n}^\pm)\/n^2}),\qq \z_{0,2n}^\pm=\l_{0,2n}^\pm-\l_{n},
\]
\[
1=\D(\l_{n}^\pm)=\D(\l_{n})+A_{n}^\pm, \qq A_{n}^\pm=\ddot \D(\l_{n}) {{\z_{n}^\pm}^2\/2}(1+{O(\z_n^\pm)\/n^2}),\qq \z_n^\pm=\l_{n}^\pm-\l_{n},
\]
which yields $\D(\l_{0,n}^\pm)-1=-A_{n}^\pm
+A_{0,n}^\pm$. Substituting \er{asDY} into the identity $\D(\l_{0,n}^\pm)-1=\sqrt{1+{\D_-^2(\l_{0,n}^\pm)\/9}}-1$ we obtain
\[
\lb{asD0n}
\D(\l_{0,n}^\pm)-1=\sqrt{1+{\D_-^2(\l_{0,n}^\pm)\/9}}-1={(q_{sn}+O(n^{-1}))^2\/18(2\pi n)^2}.
\]
Substituting \er{asD0n} and \er{asDY} into the identity
$\D(\l_{0,n}^\pm)-1=-A_{n}^\pm +A_{0,n}^\pm$ we get
$$
{\hat q_{sn}^2\/9}+{O(|\hat q_n|)\/n}={\z_{n}^\pm}^2-{\z_{0,n}^\pm}^2+{O(1)\/n^2}.
$$
Using $\l_n^\pm=(\pi n)^2+q_0\pm |\hat q_n|+O(n^{-1})$
from [MO] and $\l_n=(\pi n)^2+q_0+O(n^{-1})$ from [Ko]
we obtain $\z_n^\pm=\l_{n}^\pm-\l_{n}=\pm |\hat q_n|+O(n^{-1})$  which yields
\[
\z_{0,n}^\pm=\pm\sqrt{{\z_{n}^\pm}^2-{\hat q_{sn}^2\/9}+{O(|\hat q_n|)\/n}}=\pm
\sqrt{{|\hat q_n|^2}-{\hat q_{sn}^2\/9}}+O(n^{-1}).
\]

(iii) Let $\m_n,\n_n$ lay on the different edge of the gap
$\g_{n}$. Then Lemma \ref{T41} (i)  gives
$\g_{n}=\g_{0,2n}$.

Conversely, let $\g_{n}=\g_{0,2n}=(a,b)$ and let $\l\in\{a,b\}$. 
Then the identity $\D_0=2\D^2+{\vt'(1,\cdot)\vp(1,\cdot)\/4}-1$ yields
$1=\D_0(\l)=1+{\vt'(\l)\vp(\l)\/4}$, thus $\vt'(\l)\vp(\l)=0$,
which gives (iii)

(iv) Let $\g_n=(a,b)$ and let $a<\m_n<\n_n<b$. The proof
of other cases is similar. Let
$
\D_0=f_0+f$, where $f_0=2\D^2-1, f={\vt'(1,\l)\vp(1,\l)\/4}$.
We obtain 
$$
f_0|_{\g_n}>1,\ \ \ \ \ \ f|_{\s}>0,\ 
\ \ \ \ \ \ f|_{\g_n\sm \s}<0
\ \ \ where\ \ \ \s=(\m_n,\n_n),
$$
which yields $\D_0|_{\s}>1$ and $\D_0(\l)<1$ for any $\l\in \{a,b\}$.
Then there exist two points $\l_n^\pm\in (a,b)$ such that
$\D_0(\l_n^\pm)=1$. Thus we have $\g_{0,2n}\ss \g_n$ for all $n\ge 1$.

We will show $\g_{0,2n}=\es$ iff $\g_{n}=\es$.
If $\g_{n}=\es$, then $\g_{0,2n}\ss \g_n$ for all $n\ge 1$
yields $\g_{0,2n}=\es$.

If $\g_{0,2n}=\es$, then for $\l=\{\m_n,\n_n\}, n\ge 1$
we obtain $\D_0(\l)=2\D^2(\l)+{\vt'(1,\l)\vp(1,\l)\/4}-1=
2\D^2(\l)-1=0$, which  yields $\D^2(\l)=1$. Thus 
we deduce that $\g_{0,2n}=\es$ iff $\g_{n}=\es$.

(v) Let $q\in L^2_{even}(0,1)$. Then we get
$\vp'(1,\cdot)=\vt(1,\cdot)=\D$ (see p. 8, [23]), which
together with \er{DY} yields $\D_0={9\/4}\D^2-{5\/4}$.
Then the zeros of $\D(\l)$ and $\D_0'(\l)$(at $\D_0(\l)<0$) coincide, since all zeros of $\D_0'$ are simple.

Conversely, let $\D_0(\l_{0,n})=-{5\/4}$ for all odd $n\ge 1$.
Then Lemma \ref{T41} ii) gives $\vt(1,\l_{0,n})=\D(\l_{0,n})=\vp'(1,\l_{0,n})=0$ for all odd $n\ge 1$. Then $\e_{n}=\e_{0,2n-1}, n\ge 1$, which implies
$\vp'(1,\cdot)=\vt(1,\cdot)=\D$. 
Thus the Wronskian identity gives $\D^2-\vt'(1,\cdot)\vp(1,\cdot)=1$
and all zeros of $\D^2(\l)=1$ coincide with the Dirichlet
and Neumann eigenvalues. Then the results of [GT] or [KK]
imply $q\in L^2_{even}(0,1)$.
It is simple fact in the inverse spectral theory
and it can be proved using other methods, see [PT]. 
\BBox

\section{ The Lyapunov function $\D_k, k=1,..,m$}
\setcounter{equation}{0}

{\bf Proof of Theorem \ref{T5}}.
(i) We determine the equation for periodic eigenvalues.
Using \er{T1-2} and $4c_k^2=(1+s^k)(1+s^{-k})$ for $k\in\ol m=\{1,2,..,m\}$, we have
$$
\det (\cM_k \mp I_2)=1 \mp \Tr M_k+s^{-k}=
1+s^{-k}\mp 4{\D_0+s_k^2\/(1+s^k)}={4(c_k^2\mp (\D_0+s_k^2))\/(1+s^k)},
$$
which yields the equation  
\[
\D_0(\l)=c_k^2-s_k^2=\cos {2\pi k\/N}\in (-1,1)\sm\{0\} \ \ \ \ {\rm for \ periodic \ eigenvalues}\ \ \l,
\]
\[
\D_0(\l)=-c_k^2-s_k^2=-1 \  \ \ \ \ {\rm for \ anti-periodic \ eigenvalues}\ \ \l.
\]
Then by Lemma \ref{TLas}, all zeros of $\det (\cM_k \mp I_2)=0$
are real and simple.
  Using $\x_k=\D_0+s_k^2,\ \ \ \r_k={s_k^2\/c_k^2}(c_k^2-\x_k^2)$, we get (recall  $\dot u={\pa\/\pa \l}u$)
\[
\lb{4dD}
\D_k=\x_k+\sqrt{\r_k}=\D_0+s_k^2+\sqrt{\r_k}, \qqq
\dot\D_k=\dot\x_k+{\dot \r_k\/2\sqrt{\r_k}}=\rt(1-{s_k^2\/c_k^2}{\x_k\/\sqrt{\r_k}}\rt)\dot\D_0.
\]
For periodic eigenvalues $\l=\l_{k,n}^\pm,$  where $n\ge 0$
is even, we have $\D_k(\l)=1$ and  
\[
\x_k=\D_0+s_k^2=c_k^2-s_k^2+s_k^2=c_k^2,\ \ \ \
\sqrt{\r_k}=1-\x_k=1-c_k^2=s_k^2, \ \ at \ \l=\l_{k,n}^\pm,
\]
 and then \er{4dD} yields at $\l=\l_{k,n}^\pm$: 
$$
{\dot\D_k\/\dot\D_0}=1-{s_k^2\/c_k^2}{\x_k\/\sqrt{\r_k}}=
1-1=0,\qq
\ddot\D_k=-{s_k^2\dot\D_0\/c_k^2}\rt({\dot \D_0\/\sqrt{\r_k}}-{\x_k\dot\r_k\/2\r_k^{3/2}}\rt)=-{\dot\D_0^2\/c_k^2}\rt(1+{c_k^2\/s_k^2}\rt)<0,
$$
which gives \er{espk} for periodic eigenvalues $\l_{k,2n}^\pm$.

For anti-periodic eigenvalues $\l=\l_{k,n}^\pm,$ where $n\ge 1$
is odd, we nave $\D_k(\l)=-1$ and 
\[
\x_k=\D_0+s_k^2=-1+s_k^2=-c_k^2,\ \ \ \
\sqrt{\r_k}=-1-\x_k=-s_k^2,\ \ at \ \l=\l_{k,n}^\pm,
\ \ 
\]
and then \er{4dD} yields
$$
{\dot\D_k\/\dot\D_0}=1-{(1-c_k^2)\/c_k^2}{(1-s_k^2)\/s_k^2}=0,
\qq
\ddot\D_k=-{s_k^2\dot\D_0\/c_k^2}\rt({\dot \D_0\/\sqrt{\r_k}}-{\x_k\dot\r_k\/2\r_k^{3/2}}\rt)={\dot\D_0^2\/c_k^2}\rt(1+{c_k^2\/s_k^2}\rt)>0,
$$
which yields \er{espkI} for anti-periodic eigenvalues $\l_{k,n}^\pm$.

Let $\l=\l_{k,2n}^\pm, k\in \ol m, n\ge 0$, be periodic eigenvalues. Recall that they satisfy $\D_0(\l)=\cos {2\pi k\/N}\in (-1,1)$. Using Lemma \ref{TLas} for $c=\cos {2\pi k\/N}$ we have \er{Tas-2}.

(ii) We determine the equation for resonances for $k\in \ol m$. We have:
 $\r_k(\l)=0 \lra  \x_k(\l)=\pm c_k$. Consider the first case $+$: $n$ is {\bf even}, i.e., ${n\/2}\in \Z$. We obtain
\[
 \x_k(\l)=c_k\ \ \ \ \  \lra  \ \  \D_0(\l)>a_k^+=c_k-s_k^2\in (-1,1),\qq
  k\in \ol m.
\]
Then the resonances $r_{k,n}^\pm$  are zeros of Eq. $\D_0(\l)=c_k-s_k^2\in (-1,1),\l\in \C$. 
From Lemma \ref{TLas} we deduce that all these resonances 
are real. Moreover, we have $1>a_1^+>a_2^+>...>a_m^+$, which yields
for even all\ even $ n\ge 1$:
$$
r_{1,0}^+<r_{2,0}^+<...<r_{m,0}^+,\ \  
r_{m,n}^-...<r_{2,n}^-<r_{1,n}^-
<\l_{0,n}^-<\l_{0,n}^+<
r_{1,n}^+<r_{2,n}^+<...<r_{m,n}^+,
$$
\[
\lb{clgap}
\g_{m,n}\ss\g_{m-1,n}\ss...\ss\g_{1,n}\ss\g_{0,n},\ \ \ 
\ \ \  and\ \ G_{n}=\cap_0^m \g_{k,n}=\g_{0,n}.
\]

Consider the second case $-$: $n$ is {\bf odd}. We have
\[
 \x_k(\l)=-c_k\ \ \ \ \  \lra  \ \  \D_0(\l)=a_k^-=-c_k-s_k^2\in (-{5\/4},-1),\ \ \ {N\/3}\notin \Z,
\]
\[
 \x_k(\l)=-c_k\ \ \ \ \  \lra  \ \  \D_0(\l)=a_k^-=-c_k-s_k^2\in [-{5\/4},-1),\ \ \ {N\/3}\in \Z.
\]
Then from Theorem \ref{T3} we deduce that the zeros of $\r_k$
are zeros of Eq. $\D_0(\l)=-c_k-s_k^2\in [-{5\/4},-1),\l\in \C$. Moreover, by Lemma \ref{TLas}, they are real and simple.

If ${N\/3}\notin \Z$, then ${N\/3}=p+\ve $ for some integer $p\ge 1$
and $\ve \in (0,1)$. Then we obtain 
\[
\lb{gapevI}
\g_{0,n}\sps\g_{1,n}\sps\g_{2,n}\sps...\sps\g_{p,n},\ \ \ 
and \ \ \ \ \g_{p+1,n}\ss \g_{p+2,n}\ss ..\ss \g_{m,n},
\]
\[
\lb{gapev}
\g_{p,n}\ss \g_{p+1,n} \ \ or \ \ \ \ \g_{p,n}\sps \g_{p+1,n}
\ \ and \ \ G_n=\cap_1^N\g_{k,n}=\g_{p,n}\cap \g_{p+1,n},
\]
for all odd $n\ge 1$. If $p={N\/3}\in \Z$, then we get for all odd $n\ge 1$:
\[
\lb{gapod1}
\g_{0,n}\sps\g_{1,n}\sps\g_{2,n}\sps...\sps\g_{p,n},\ \ \ 
and \ \ \ \ \g_{p,n}\ss \g_{p+1,n}\ss ..\ss \g_{m,n},
\]
\[
\lb{gapod}
G_n=\cap_{k=0}^m\g_{k,n}=\g_{p,n},\ \ \D_0(r_{p,n}^\pm)=-{5\/4}. 
\]
Resonances $r_{k,2n}^\pm , k\in \ol m,n\ge 0 $ satisfies 
$\D_0(r_{k,2n}^\pm )=c_k-s_k^2\in (-1,1)$.
Using Lemma \ref{TLas} for $c=c_k-s_k^2$ we obtain
asymptotics \er{Tas-4} for $r_{k,2n}^\pm $.

Resonances $r_{k,2n+1}^\pm , k\in \ol m, k\ne {N\/3},n\ge 0$ satisfies 
$\D_0(r_{k,2n+1}^\pm )=-c_k-s_k^2\in (-{5\/4},-1)$.
Using Lemma \ref{TLas} for $c=c_k-s_k^2$ we obtain
asymptotics \er{Tas-4} for $r_{k,2n+1}^\pm $.

(iii) Finally, we obtain the following identities   for $(n,k)\in \N\ts \ol m$: 
\[
S_n=\cup_{k=0}^m\s_{k,n},\ \ \s_{0,n}=(\l_{0,n-1}^+,\l_{0,n}^-),\ \ \ \s_{k,n}=(r_{k,n-1}^+,r_{k,n}^-),\ \ \ \ \ 
\]
\[
G_n=\cap_{k=0}^m\g_{k,n},\  \g_{0,n}=(\l_{0,n}^-,\l_{0,n}^+),\ \ \ \ \g_{k,n}=(r_{k,n}^-,r_{k,n}^+),
\]
which together with \er{clgap}, \er{gapevI}-\er{gapod}  yield gives \er{T5-1}-\er{T5-3}.

(iv) Using Theorem \ref{T4} we deduce that
$G_{2n}=\g_{0,2n}=\es$ iff $\g_{n}=\es $. Moreover, 
asymptotics $|\g_{n}|\to 0$ as $n\to \iy$  give
$|G_{2n}|\to 0$ as $n\to \iy$.

(v) In the case $p={N\/3}\in \Z$ and odd $n\ge 1$
using the identity \er{gapod},
we have $G_n=\g_{p,n}=(r_{p,n}^-,r_{p,n}^+)$,
where $r_{p,n}^\pm$ are zeros of Eq. $\D_0(\l)=-{5\/4}$.
Theorem \ref{T4} (v) gives $r_{p,n}^\pm=\l_n$  iff
$q\in L_{even}^2(0,1)$.

We determine the asymptotics \er{Tas-6}.
Let $(z_n^\pm)^2=r_{p,n}^\pm$. Using ${9\D^2-\D_-^2-5\/4}={-5\/4}$ we have
$9\D^2(r_{p,n}^\pm)=\D_-^2(r_{p,n}^\pm)$. Thus in the case $q=0$ we have $\cos z=0$ and the
corresponding zeros are given by $z_n=\pi (n-{1\/2}),n\ge 1$. 

 In the case $q\ne 0$ we have $\D(\l)=\pm {1\/3}\D_-(\l)={o(1)\/\sqrt \l}$
 as $\l\to \iy$.  Then by Lemma \ref{T41},\
$|z_n^{\pm}-z_n|<{\pi\/4}$ for all $n>n_0,$ for some $n_0\ge 1$.
  Thus we obtain  $z_{n}^{\pm}=z_n+\z_n^\pm,\z_n^\pm\to 0$.
 Moreover, using \er{asLq}, we get ${o(1)\/n}=\D(z_{n}+\z_n^\pm)=\cos (z_{n}+\z_n^\pm-{q_0\/2z_n}+{o(1)\/n^2})$,
 which implies $\z_n^\pm={q_0\/2z_n}+{o(1)\/n}$. 
Let $\z_n^\pm={q_0+v\/2z_n}, v_n^\pm\to 0$
and for $R_n^\pm={v_n^\pm\/2z_n}+{o(1)\/n^2}$
using \er{asLq} we obtain
$$
\D(z_{n}+\z_n^\pm)=\cos (z_{n}+R_n^\pm)
=(-1)^{n-1}\cos (-{\pi\/2}+R_n^\pm)=(-1)^{n-1}\sin ({v_n^\pm\/2z_n}+{o(1)\/n^2})
$$
 Using this and \er{asY} we get 
 \[
 \D_-(r_{p,n}^\pm)=(-1)^n{\wt q_{cn}\/2z_n}+{O(1)\/n^2},\ \ \qq  n\to \iy.
\]
Thus the Eq. $\D(r_{p,n}^\pm)=\pm {1\/3}\D_-(r_{p,n}^\pm)$
yields \er{Tas-6}.

(vi) In the case ${N\/3}=p+\ve $, where $p\ge 1$
is integer and $\ve \in (0,1)$ and odd $n\ge 1$
the identity \er{gapev} we have $G_n=\g_{p,n}\cap \g_{p+1,n}.$
In \er{Tas-4} we will show that $|\g_{k,n}|\to \iy$ as $n\to \iy, k\ne {N\/3}$, which yields $|G_n|\to \iy$  as $n\to \iy$.
\BBox

We will prove the last Corollary \ref{T6}.

{\bf Proof of Corollary \ref{T6}}. 
(i) Recall that by Theorem \ref{T5}, each $G_{2n}=\g_{0,2n}\ss\g_n, n\ge 1$.
Thus by Theorem \ref{T5} (iv), $G_{2n}=\g_{0,2n}=\es$ for all $n\ge n_0$ and for some $n_0\ge 1$ iff $q$ is a finite gap potential for the operator
$-y''+qy$ on the real line.

Consider the odd gaps $G_{2n-1}, n\ge 1$.
Recall that by Theorem \ref{T5}, if ${N\/3}\notin \N$, then
$|G_{2n-1}|\to \iy$ as $n\to \iy$. 

Assume that $p={N\/3}\in \N$. In this case by Theorem \ref{T5},
 $G_{2n+1}=\g_{p,2n+1}=\es$ iff $q\in L_{even}^2(0,1)$.
 Thus the statement i) has been proved.

(ii) The function $f_\o$ satisfies the Eq.
$-f_\o''+qf_\o=E_nf_\o$ on the interval $[0,1]$
and $f_\o(0)=f_\o(1)=0, f_\o'(1)^2=e^{2h_n}>0$
for some $h_n\in \R$. Each constant $h_n$ is so-called
norming constant for the Sturm-Liouville problem
$-y''+qy=\l y$ on the interval $[0,1]$
with the Dirichlet boundary conditions
$y(0)=y(1)=0$ [PT]. 

Recall that $\m_n, n\ge 1,$ is the Dirichlet spectrum of 
the problem $-y''+qy=\l y, y(0)=y(1)=0$ on the unit interval
$[0,1]$ and $\m_n=E_n, n\ge 1$.

Recall the well know result from [PT]: the mapping 
$
\gF:q\to \lt(q_0,(\vk_n(q))_{1}^\iy;(h_n(q))_{1}^\iy\rt)
$
is a real-analytic isomorphism between $L^2(0,1)$ and $\R\ts\cK\ts\ell^2_1$. This gives the statement ii).

(iii) If potential $q$ is even, i.e., $q\in L_{even}^2(0,1)$,
then each $h_n=0, n\ge 1$ (see [PT]). 

Recall the well know result from [PT]: the mapping 
$
\gF_e:q\to \lt(q_0,(\vk_n(q))_{1}^\iy\rt)
$
is a real-analytic isomorphism between $L_{even}^2(0,1)$ and $\R\ts\cK$. This gives the statement iii).
\BBox

 \no {\bf Acknowledgments.}
E. Korotyaev was partly supported by DFG project BR691/23-1.
The various parts of this paper were written at the Mittag-Leffler Institute, Stockholm  and in the Erwin Schr\"odinger Institute for Mathematical Physics, Vienna, the first author is grateful to the Institutes for the hospitality. The first author would like to thank Markus Klein for useful discussions.

\end{document}